\documentclass[12pt,a4paper,draft]{article}

\usepackage{amsmath,amsfonts,amsthm,amssymb}
\usepackage{mathrsfs}

\usepackage{a4wide}
\newtheorem{theorem}{Theorem} 

\newtheorem{remark}{Remark}
\newtheorem{example}{Example}
\newtheorem{proposition}{Proposition}

\let\scr\mathscr

\def\Pb{\mathbf{P}}

\def\Ex{\mathbf{E}}      
\def\AA{\mathbb{A}}

\def\KK{\mathbb{K}}

\def\UU{\mathbb{U}}

\def\1{\mbox{1\hspace{-.25em}I}}
\def \Liminf{\mathop{\underline{\lim}}\limits}

\begin{document}

\title{Volatility Estimation of Hidden Markov Processes and Adaptive Filtration}
\author{Yury A. Kutoyants\\
{\small Le Mans  University,  France}, \\
{\small Tomsk State University, Tomsk, Russia}}


\date{}

\maketitle
\begin{abstract}
The partially observed linear Gaussian system of stochastic differential
equations with low noise in observations is considered. The coefficients of
this system are supposed to depend on some unknown parameter.  The problem of
estimation of these parameters is considered and the possibility of the
approximation of the filtering equations is discussed. An estimators are used
for estimation of the quadratic variation of the derivative of the limit of
the observed process. Then this estimator is used for nonparametric estimation
of the integral of the square of volatility of unobservable component. This
estimator is also used for construction of method of moments estimators in the case
where the drift in observable component and the volatility of the state
component depend on some unknown parameter.  Then this method of moments estimator
and Fisher-score device allow   us to introduce the One-step MLE-process and
adaptive Kalman-Bucy filter. The asymptotic efficiency of the proposed filter
is discussed.

\end{abstract}


\section{Introduction}

Partially observed systems and especially different Kalman filter models
occupy an important place in statistics due to a large diversity of applied
problems where such statistical models are used. There is extensive engineering
literature devoted to identification problems for such partially observed
discrete time 
models. Note that the continuous time systems are less studied. 

This work is devoted to the construction of adaptive Kalman-Bucy filter for
partially observed linear system with the small noise in observation
equation. It is supposed that the coefficients of these equations depend on
some unknown parameter and the adaptive filter is constructed in several
steps. First we propose a nonparametric estimator of some integral functional,
then this estimator is used for construction of method of moments estimator of the
unknown parameter. This estimator is used as preliminary for construction of
One-step MLE-process and finally this process is used for writing the adaptive
Kalman-Bucy filter. The problems of parameter estimation by observations of
partially observed linear systems with small noise in observations and in
state equation for more complex (multidimensional, nonlinear) models were
extensively studied in filtration theory, where approximate filters were
proposed and studied in different asymptotics, see, e.g.,
\cite{CMT05},\cite{EAM95},\cite{EM02} and references there in. In particular
the order of the errors of approximations are obtained for the large diversity
of the statement of the problems, see \cite{FP89}, \cite{KBS84}, \cite{Pic91},
\cite{Pic93}. The similar problems with hidden telegraph process were studied
in \cite{PCh09} and \cite{KhK18}.  The statistical problems for partially
observed linear and non linear systems were studied in \cite{Kut94}, Chapter 6
and in \cite{SKh96}.  Note that in \cite{Kut94} and in \cite{SKh96} it is
supposed that the small noises are in the observation and state equations.
The parameter estimation problems in the case of small noise in observations
only were considered recently in \cite{Kut19}. The construction of
asymptotically optimal estimator of the parameter in volatility of unobserved
component for partially observed system by discrete time observations was
studied in \cite{GJ1-01},\cite{GJ2-01}.  The asymptotic behavior of the filter
for such partially observed nonlinear continuous time system with small noise
in observations was studied in the mentioned above works \cite{Pic91},
\cite{Pic93}. 

Let us give more detailed exposition of the presented in this work study. We
consider the linear two-dimensional partially observed system
\begin{align}
\label{22-41}
{\rm d}Y_t&=a\left(\vartheta ,t\right)Y_t{\rm d}t+b\left(\vartheta ,t\right){\rm d}V_t,\quad
Y_0=y_0,\quad 0\leq t\leq T,\\
\label{22-42}
{\rm d}X_t&=f\left(\vartheta ,t\right)Y_t{\rm d}t+\varepsilon \sigma
\left(t\right){\rm d}W_t,\qquad X_0=0,\quad 0\leq t\leq T,
\end{align}
where the Wiener processes $V_t,0\leq t\leq T$ and $W_t,0\leq t\leq T$ are
 independent. The solution $Y^T=\left(Y_t,0\leq t\leq T\right)$
of the state equation \eqref{22-41} can not be observed directly and we have
available the observations $X^T=\left(X_t,0\leq t\leq T\right)$ only.  Here
$a\left(\cdot \right), b\left( \cdot \right),$ $ f\left(\cdot \right)$ and
$\sigma \left(\cdot \right)$ are some bounded functions and $\varepsilon \in
(0,1]$ is {\it small} parameter. 

The conditional expectation $m\left(\vartheta ,t\right)=\Ex_\vartheta
\left(Y_t|X_s,0\leq s\leq t\right)$ satisfies the equations of Kalman-Bucy
filtration \cite{KB61}, \cite{LS01}
\begin{align}
\label{22-43}
{\rm d}m\left(\vartheta ,t\right)&=a\left(\vartheta ,t\right)m\left(\vartheta
,t\right){\rm d}t +\frac{\gamma \left(\vartheta ,t\right)f\left(\vartheta
  ,t\right)}{\varepsilon ^2\sigma \left(t\right)^2}\left[{\rm
    d}X_t-f\left(\vartheta ,t\right)m\left(\vartheta ,t\right){\rm d}t\right]
,\\
\frac{\partial \gamma \left(\vartheta ,t\right)}{\partial t}&=2a
\left(\vartheta ,t\right)\gamma \left(\vartheta ,t\right)-\frac{\gamma
  \left(\vartheta ,t\right)^2f \left(\vartheta ,t\right)^2 }{\varepsilon
  ^2\sigma \left(t\right)^2} +b\left(\vartheta ,t\right)^2,\qquad \gamma
\left(\vartheta ,0\right)=0,
\label{22-44}
\end{align}
with the initial value $m\left(\vartheta ,0\right)=y_0 $. If the value of
$\vartheta $ is unknown, then to approximate $m\left(\vartheta ,t\right)$ we
can use some estimator $\bar\vartheta _\varepsilon $ and to use something like
$m\left(\bar\vartheta _\varepsilon,t\right) $, but there are some
problems. For example,  if
$\bar\vartheta _\varepsilon $ is the MLE $\hat\vartheta _\varepsilon$
constructed by observations $X^T$, then there is a problem of definition of
the It\^{o} integral in \eqref{22-43}
\begin{align*}
\int_{0}^{t}\frac{\gamma (\hat\vartheta
  _\varepsilon,s)f(\hat\vartheta _\varepsilon
  ,s)}{\varepsilon ^2\sigma \left(s\right)^2}{\rm d}X_s
\end{align*}
because the estimator $\hat\vartheta _\varepsilon $ depends on the whole
trajectory $X^T$. The using a series of MLEs $\hat\vartheta _{t,\varepsilon }$
constructed by observations $X^t=\left(X_s,0\leq s\leq t\right)$ for each $\in
(0,T]$ can provide a good approximation 
 of $m\left(\vartheta ,t\right)$ but the numerical realization of such
  algorithm can be a difficult problem. The goal of this work is the
  construction of another estimator $\vartheta _{t,\varepsilon }^\star,0<t\leq T $,
  called One-step MLE-process, which depends on the observations $X^t $, can
  be easily calculated and which is asymptotically ($\varepsilon \rightarrow
  0$) equivalent to the MLE $\hat\vartheta _{t,\varepsilon },0<t\leq T $. Then this estimator can be used for approximation
  $m\left(\vartheta ,t\right)$ as follows $m\left(\vartheta _{t,\varepsilon
  }^\star ,t\right)$. We describe the error of approximation $m\left(\vartheta
  ,t\right)-m\left(\vartheta _{t,\varepsilon }^\star ,t\right) $ and discuss
  the optimality of such approximation.
  
The construction of $\vartheta _{t,\varepsilon }^\star$ requires a
preliminary consistent estimator $\bar\vartheta_{\tau ,\varepsilon } $, which
we define with the help of nonparametric estimator of the quadratic variation
of the trend coefficient of the observed process. Introduce the function
\begin{align*}
\Psi _\tau\left(\vartheta \right)=\int_{0}^{\tau}f\left(\vartheta
,s\right)^2b\left(\vartheta ,s\right)^2{\rm d}s,\qquad \vartheta \in\Theta
,\quad \tau \in (0,T]. 
\end{align*}

 Note that we have formally 
\begin{align*}
\left.\frac{\partial X_t}{\partial t}\right|_{\varepsilon =0}=f\left(\vartheta
,t\right)Y_t .
\end{align*}
 Let us denote the trend coefficient of $X_t$ as
 $M_t=f\left(\vartheta ,t\right)Y_t$ and put  $S\left(\vartheta
 ,t\right)=f\left(\vartheta ,t\right)b\left(\vartheta ,t\right) $. The process
  $M_t$ has the stochastic differential
\begin{align*}
{\rm d}M_t=\left[f'\left(\vartheta ,t\right)-a\left(\vartheta ,t\right)f\left(\vartheta ,t\right)\right]Y_t{\rm
  d}t+S\left(\vartheta ,t\right){\rm d}V_t, \qquad M_0=f\left(\vartheta ,0\right)y_0
\end{align*}
and therefore    $M_\tau^2$ admits  the representation 
\begin{align*}
M_\tau^2=2\int_{0}^{\tau}M_t\,{\rm
  d}M_t+\int_{0}^{\tau}S\left(\vartheta ,t\right)^2{\rm d}t. 
\end{align*}
 The quadratic variation of  $M_t$ can be calculated as follows
\begin{align*}
\Psi _\tau\left(\vartheta \right)=M_\tau^2-2\int_{0}^{\tau}M_s\,{\rm
  d}M_s,\qquad 0\leq \tau\leq T. 
\end{align*}
The preliminary estimator $\bar\vartheta_{\tau ,\varepsilon }  $ is
constructed in two steps. First we estimate the derivative of $X_t$, then we
estimate the quadratic variation of this derivative. Introduced statistic
$\Psi _{\tau ,\varepsilon } $ with a fixed $\tau $ allows us to realize these
two procedures simultaneously. Using this statistic and relation $\Psi _{\tau
  ,\varepsilon }=\Psi _\tau\left(\bar\vartheta_{\tau ,\varepsilon } \right)$
we obtain the method of moments estimator $\bar\vartheta_{\tau ,\varepsilon }  $.
Of course, the process $X_t$ is not differentiable and we propose, as usual in
nonparametric statistics, a {\it slow differentiation}.

 The asymptotic properties of the maximum 
likelihood estimators (MLE) and Bayesian estimator (BE) of this parameter for 
the model \eqref{22-41}-\eqref{22-42}, were 
described in the work \cite{Kut19}. Recall that  the MLE $\hat\vartheta
_\varepsilon $ and BE $\tilde\vartheta _\varepsilon $ are  defined by the relation
\begin{align}
\label{22-45a}
L(\hat\vartheta
_\varepsilon,X^T)=\sup_{\vartheta \in\Theta }L\left(\vartheta
,X^T\right),\qquad \quad \tilde\vartheta _\varepsilon =\frac{ \int_{\Theta
  }^{}\vartheta p\left(\vartheta  \right)L\left(\vartheta ,X^T\right){\rm
    d}\vartheta }{\int_{\Theta }^{} p\left(\vartheta  \right)L\left(\vartheta
  ,X^T\right){\rm d}\vartheta }, 
\end{align}
where the likelihood ratio is (see \cite{LS01})
\begin{align}
\label{22-45}
L\left(\vartheta ,X^T\right)=\exp\left(\int_{0}^{T}\frac{f\left(\vartheta
  ,t\right)m\left(\vartheta ,t\right)}{\varepsilon ^2\sigma
  \left(t\right)^2}\;{\rm d}X_t-\int_{0}^{T}\frac{f\left(\vartheta
  ,t\right)^2m\left(\vartheta ,t\right)^2}{2\,\varepsilon ^2\sigma
  \left(t\right)^2}\;{\rm d}t\right).
\end{align}

It was shown that these estimators are
consistent, asymptotically normal 
\begin{align}
\label{22-46}
\frac{\hat\vartheta _\varepsilon -\vartheta }{\sqrt{\varepsilon
}}\Longrightarrow  {\cal N}\left(0,{\rm I}\left(\vartheta
\right)^{-1}\right),\qquad \quad {\rm I}\left(\vartheta \right)=
\int_{0}^{T}\frac{\dot S\left(\vartheta ,t\right)^2}{2S\left(\vartheta
  ,t\right)\sigma \left(t\right)}\;{\rm d}t ,
\end{align}
the polynomial moments converge and the both estimators are asymptotically
efficient. Here in the sequel   dot means derivation w.r.t. $\vartheta $ and
prim means derivation w.r.t. $t$.

 The model \eqref{22-41}-\eqref{22-42} has some unusual features. For example,
 the {\it empirical} Fisher information
\begin{align*}
{\rm I}_\varepsilon \left(\vartheta \right)=\int_{0}^{T}{\left[\frac{\dot f\left(\vartheta ,t\right)m\left(\vartheta
    ,t\right)+f\left(\vartheta ,t\right)\dot m\left(\vartheta
    ,t\right)}{\sigma \left(t\right)  }\right]^2}\;{\rm d}t \longrightarrow
0,\qquad {\rm as}\quad \varepsilon \rightarrow 0
\end{align*}
 and for all $t\in (0,T]$
  there is  convergence in distribution
\begin{align*}
\varepsilon ^{-1/2}\left[\dot f\left(\vartheta ,t\right)m\left(\vartheta
  ,t\right)+f\left(\vartheta ,t\right)\dot m\left(\vartheta
  ,t\right)\right]\Longrightarrow {\frac{\dot
    S\left(\vartheta ,t\right)  \sqrt{\sigma
    \left(t\right)}}{\sqrt{2S\left(\vartheta ,t\right)}}}\;  \xi _t.
\end{align*}
Here  $\xi _t, t\in (0,T]$ are mutually independent standard
  Gaussian (${\cal     N}\left(0,1\right)$) random variables. This convergence
  provides convergence in probability 
\begin{align*}
\varepsilon ^{-1}{\rm I}_\varepsilon \left(\vartheta
\right)=\int_{0}^{T}{\left[\frac{\dot f\left(\vartheta
      ,t\right)m\left(\vartheta ,t\right)+f\left(\vartheta ,t\right)\dot
      m\left(\vartheta ,t\right)}{\sqrt{\varepsilon }\sigma \left(t\right)
    }\right]^2}\;{\rm d}t \longrightarrow {\rm I}\left(\vartheta \right).
\end{align*}

To calculate just  one value of $L\left(\vartheta ,X^T\right) $ by \eqref{22-45} we
need the solutions $m\left(\vartheta ,t\right),\gamma \left(\vartheta
,t\right),$ $0\leq t\leq T$ of the equations \eqref{22-43},\eqref{22-44} and
therefore, as we wrote above, the numerical calculation of the MLE and BE can be
difficult problem.

 Here we propose a different construction of the  estimator  of  this parameter
 which is computationally simpler than the  MLE and which is asymptotically
 equivalent to the MLE \cite{Kut19}.

\section{Model of  observations}
Consider a {non-homogeneous} partially observed linear system described by the
equations
\begin{align}
\label{23-1}
{\rm d}X_t &= f\left(\vartheta,t\right)\,Y_t\,{\rm d}t+\varepsilon
\sigma\left(t\right)\, {\rm d}W_t,\quad \qquad X_0=0, \\
 {\rm d} Y_t &=
a\left(\vartheta,t\right)Y_t\,{\rm d}t+b\left(\vartheta,t\right)\, {\rm d}V_t,
\qquad Y_0= y_0,
\label{23-2}
\end{align}
 where $f\left(\cdot \right), \sigma \left(\cdot \right), a\left(\cdot
 \right)$ and $b\left(\cdot \right)$ are known, smooth functions, while $W_t,
 0\leq t\leq T$ and $V_t, 0\leq t\leq T$ are two independent Wiener
 processes. The process $X^T=\left(X_t,0\leq t\leq T\right)$ is observed and
 the process $Y^T=\left(Y_t,0\leq t\leq T\right)$ is {\it hidden}. We are
 interested by the same problems: estimation of the finite-dimensional parameter
 $\vartheta\in\Theta \subset {\cal R}^d  $ and the construction of adaptive
 filter. 

Here we study the method of moments estimators (MME) based on some
nonparametric estimation of one integral, One-step MLE-process and a adaptive filter
based on this estimator-process.

The limit model ($\varepsilon =0$) corresponds to the observations of
the integral of the hidden process and the possibility of the consistent
estimation is equivalent to the possibility to estimate without error the
parameter of the observed Gaussian process.

\section{Limit model}
Let us now present how an estimator $\vartheta _\varepsilon ^*$ of
$\vartheta_0 $ (true value) could be obtained without error by means of the
observations $X^T=x^T=\left(x_t,0\leq t\leq T\right)$, when $\varepsilon
=0$. We
suppose that we observe a Gaussian process $x^T$ satisfying
\begin{align*}
\frac{\partial x_t}{\partial t}&=f\left( \vartheta,t\right)Y_t, \quad
x_0=0,\qquad 0\leq t\leq T\\
 {\rm d}Y_t&=a\left( \vartheta,t\right)Y_t{\rm
  d}t+b\left( \vartheta,t\right){\rm d}V_t, \quad Y_0=y_0
\end{align*}
and aim to estimate $\vartheta_0 $. Suppose that the function
$f\left(\vartheta,t\right)$ has a continuous derivative  {with respect to}  
$t$.  Then the process $x^T$ is described by the following stochastic
differential equation
\begin{align*}
{\rm d}x_t'&= f'\left( \vartheta_0,t\right)Y_t{\rm d}t+f\left(
\vartheta_0,t\right){\rm d}Y_t\\ &= \left[f'\left( \vartheta_0,t\right)a\left(
  \vartheta_0,t\right)+f\left( \vartheta_0, t\right) \right] Y_t{\rm d}t+f\left(
\vartheta_0,t\right)b\left( \vartheta_0,t\right){\rm d}V_t,
\end{align*}
where $f'\left(\vartheta_0,t\right)$ is the derivative of $f\left(\vartheta_0
,t\right)$  {with respect to}  $t$. 

Since the volatility function depends on the unknown parameter $\vartheta_0 $
through the product 
$S \left(\vartheta_0,t\right)=f\left(  \vartheta_0,t\right)b\left(
\vartheta_0,t\right)$, it is clear that consistent estimation is
feasible. It will be shown that the identifiability condition and the Fisher information  have to be 
based mainly on the function $S \left(\vartheta_0,t\right)$. For example if
$f\left(  \vartheta,t\right)=\vartheta  $  and $b\left(
\vartheta,t\right)={1}/{\vartheta }$, then consistent estimation is
impossible and the (limit) Fisher information is equal to $0$.  

Recall that by It\^o's formula, we have the following relation for the
quadratic variation 
\begin{align}
\label{23-MME-1}
{x_t'}^2=2\int_{0}^{t}x_s'{\rm d}x_s+\int_{0}^{t}b\left(
\vartheta_0,s\right)^2 f\left( \vartheta_0, s\right)^2{\rm d}s
\end{align}
and the function
\begin{align*}
\Psi_t ={x_t'}^2-2\int_{0}^{t}x_s'{\rm
  d}x_s'=\int_{0}^{t}b\left(  \vartheta_0,s\right)^2 f\left( 
\vartheta_0, s\right)^2{\rm d}s\equiv K\left(\vartheta_0,t\right)
\end{align*}
is deterministic. Under mild
identifiability conditions, the observed function $\Psi_t$
defines $\vartheta_0 $ for any $t\in (0,T]$ without 
  error. For example, the estimator $\vartheta ^*$ defined by the equation
$K\left(\vartheta^* ,t\right)=\Psi_t$
is without error, i.e., $ \vartheta^*=\vartheta_0 $.

For example, if $f\left(\vartheta ,t\right)=f\left(t\right)$,
$b\left(\vartheta ,t\right)^2=h\left(t\right)+\vartheta g\left(t\right) $,
where all functions and $\vartheta $ are positive, then the ``estimator'' 
\begin{align*}
\vartheta_t^*=\left(\int_{0}^{t}g\left(s\right)f\left(s\right)^2{\rm
  d}s\right)^{-1}    \left[\Psi_t-  \int_{0}^{t}   h\left(s\right)f\left(s\right)^2{\rm d}s\right]=\vartheta _0.
\end{align*}

Therefore the consistent estimation by observations \eqref{23-1} is not
excluded.

{Note} that if we consider the limit
  model, then the measures that correspond to observations with
  different $\vartheta $ are
  singular. Moreover $\Psi_t$ does not depend on
  $a\left(\vartheta_0,t\right)$.

\section{Consistent estimator}

One way to estimate $\vartheta $ using observations from \eqref{23-1} is by
estimating first the function $\Psi_t $ and then using the estimator, $\Psi
_{t,\varepsilon}$, to estimate $\vartheta $, i.e., to solve for some $t$ the
equation
\begin{align}
\label{23-MME-2}
\Psi _{t,\varepsilon}=\int_{0}^{t}b (  \check\vartheta_\varepsilon ,s )^2 f ( 
\check \vartheta_\varepsilon , s )^2{\rm d}s,
\end{align}
where $ \check\vartheta_\varepsilon  $ is the MME.  
Of course,   to verify the consistency of $\check\vartheta_\varepsilon $ we need
the condition of identifiability.

 {Another} possibility is to define the estimator as 
\begin{align}
\label{23-MME-3}
 \check \vartheta _{t,\varepsilon}=\arg\inf_{\vartheta \in\Theta }
\int_{0}^{t}\left[\Psi _{s,\varepsilon}- K\left(\vartheta
  ,s\right) \right]^2{\rm d}s.
\end{align}
For example, let $f\left(\vartheta,t\right)=f\left(t\right)$ and
$b\left(\vartheta,t\right)=\vartheta g \left(t\right)$. Then
\begin{align*}
\check\vartheta _{t,\varepsilon}  = \left\{  {\int_{0}^{t}\Psi
  _{s,\varepsilon}g\left(s\right){\rm d}s} {\Bigg /}
                {\int_{0}^{t}g\left(s\right)^2{\rm 
    d}s}  \right\} ^{1/2}
\end{align*} 
and from the consistency of $\Psi _{s,\varepsilon},0\leq s\leq t $ we obtain
the consistency of $\check\vartheta _{t,\varepsilon} $.

Therefore to estimate $\vartheta $ we first estimate the function $\Psi _t$
and then using this estimator of $\Psi _t$ we estimate the parameter
$\vartheta $. Below we realize this program.

\section{Nonparametric estimation.}

Suppose that we have the partially observed system
\begin{align}
\label{23_NE-6}
{\rm d}X_t&=f\left(t\right)Y_t{\rm d}t+\varepsilon \sigma
\left(t\right){\rm d}W_t,\qquad X_0=0,\qquad 0\leq t\leq T,\\
\label{23_NE-7}
{\rm d}Y_t&=a\left(t\right)Y_t{\rm d}t+b\left(t\right){\rm d}V_t,\quad \qquad
Y_0=y_0,
\end{align}
where $a\left(\cdot \right),b\left(\cdot \right),f\left(\cdot \right)$ and
$\sigma \left(\cdot \right)$ are unknown functions.

Consider the problem of estimation of the  function
\begin{align*}
\Psi _\tau =\int_{0}^{\tau }f\left(t\right)^2b\left(t\right)^2{\rm d}t,\qquad
0<\tau \leq T 
\end{align*}
by observations $X^T$. Remind that this is quadratic variation of the
derivative of the process $X_\tau $ at the  point $\varepsilon =0$. 

Introduce the statistic
\begin{align*}
\Psi _{\tau ,\varepsilon} =\sum_{i=0}^{N_{\tau,\varepsilon
}-1}\left(\frac{X_{t_{i+1}+\delta _\varepsilon }-X_{t_{i+1}}}{\delta
  _\varepsilon }-\frac{X_{t_{i}+\delta _\varepsilon }-X_{t_{i}}}{\delta
  _\varepsilon }\right)^2,\qquad\quad  0<\tau\leq T.  
\end{align*}
Here $t_i=i\varphi _\varepsilon , N_{\tau,\varepsilon
}=\left[\frac{\tau}{\varphi _\varepsilon }\right]$, the rates
$\varphi_\varepsilon \rightarrow 0,\delta _\varepsilon \rightarrow 0 $ will be
defined later. Just note that as the first step is derivation we wait that the
rate $\delta _\varepsilon \rightarrow 0$ has to be faster than the step of
discretization $\varphi
_\varepsilon \rightarrow 0$.

Let us explain why this statistic can be a consistent estimator of $\Psi_\tau $. We have
\begin{align*}
\frac{X_{t_{i+1}+\delta _\varepsilon }-X_{t_{i+1}}}{\delta _\varepsilon
}=\frac{1}{\delta _\varepsilon }\int_{t_{i+1}}^{t_{i+1}+\delta _\varepsilon
}f\left(t \right)Y_t\;{\rm d}t+\frac{\varepsilon }{\delta
  _\varepsilon }\int_{t_{i+1}}^{t_{i+1}+\delta _\varepsilon }\sigma\left(t\right)\;{\rm
  d}W_t. 
\end{align*}
Hence if we take $\varepsilon \delta _\varepsilon ^{-1/2}\rightarrow 0$, then 
\begin{align*}
\frac{X_{t_{i+1}+\delta _\varepsilon }-X_{t_{i+1}}}{\delta _\varepsilon
}= f\left( {t_{i+1}} \right)Y_{t_{i+1}}+o\left(1\right),\qquad 
\frac{X_{t_{i}+\delta _\varepsilon }-X_{t_{i}}}{\delta _\varepsilon 
}= f\left(  t_{i} \right)Y_{t_{i}}+o\left(1\right).
\end{align*}
 Further, using the smoothness of the functions we can write  formally (!)
\begin{align*}
\Psi _{\tau ,\varepsilon} &= \sum_{i=0}^{N_{\tau,\varepsilon
}-1} \left(f\left( {t_{i+1}} \right)Y_{t_{i+1}}-f\left( t_{i}\right)Y_{t_{i}}
\right)^2+o\left(1\right)= \sum_{i=0}^{N_{\tau,\varepsilon 
}-1}f\left(  t_{i} \right)^2 \left(Y_{t_{i+1}}-Y_{t_{i}} \right)^2+o\left(1\right)\\
&=\sum_{i=0}^{N_{t,\varepsilon 
}-1}f\left(  t_{i} \right)^2
\left(\int_{t_i}^{t_{i+1}}a \left( s \right)Y_s{\rm
  d}s+\int_{t_i}^{t_{i+1}}b \left( s \right){\rm d}V_s  \right)^2+o\left(1\right)\\ 
&=\sum_{i=0}^{N_{\tau,\varepsilon 
}-1}f\left( {t_{i}}\right)^2
\left(\int_{t_i}^{t_{i+1}}b\left( s \right){\rm d}V_s
\right)^2+o\left(1\right)=\sum_{i=0}^{N_{\tau,\varepsilon 
}-1}f\left( {t_{i}}\right)^2b\left( t_{i} \right)^2
\left(V_{t_{i+1}}-V_{t_i}
\right)^2+o\left(1\right)\\
&=\sum_{i=0}^{N_{\tau,\varepsilon  
}-1}f\left(  t_{i} \right)^2
b\left( {t_{i}} \right)^2  \;\left(t_{i+1}-t_{i}\right) +o\left(1\right)\\
&\longrightarrow \int_{0}^{\tau}f \left(s \right)^2
b \left(s \right)^2{\rm d}s=\Psi _\tau .
\end{align*}
Therefore this statistic can be a consistent estimator of the integral $\Psi_\tau$. 

\begin{proposition}
\label{P23-1} Let  the functions $a\left(\cdot
  \right),b\left(\cdot \right),f\left(\cdot \right), \sigma \left(\cdot
  \right)\in {\cal C}^1\left[0,T\right]$ and $\delta _\varepsilon =\varepsilon
  ,\varphi _\varepsilon =\varepsilon ^{1/3}$. Then for any $p>0$ there exists
  a constant $C>0$ such that
\begin{align}
\label{23_NE-8}
\Ex \left|\Psi _{\tau ,\varepsilon}-\Psi _{\tau  } \right|^p\leq C\tau
^{p/2}\varepsilon ^{p/3}. 
\end{align}

\end{proposition}
\begin{proof}
Now we repeat the given above expansions with description of the corresponding
errors. Throughout the proof, $C$ denote ``generic'' strictly positive
constants, which can vary from formula to formula (and even in the same
formula).  We have
\begin{align*}
\frac{X_{t_{i+1}+\delta _\varepsilon }-X_{t_{i+1}}}{\delta _\varepsilon
}&=f\left({t_{i+1}}\right)Y_{t_{i+1}}+\frac{\varepsilon }{\delta
  _\varepsilon }\int_{t_{i+1}}^{t_{i+1}+\delta _\varepsilon }\sigma\left(s\right)\;{\rm
  d}W_s\\&\qquad +\frac{1}{\delta _\varepsilon }\int_{t_{i+1}}^{t_{i+1}+\delta _\varepsilon
}\left[f\left( s \right)Y_s-f\left( {t_{i+1}}
  \right)Y_{t_{i+1}}\right]\;{\rm d}s \\
&= f\left( {t_{i}}
\right)Y_{t_{i+1}}+\left[f\left( {t_{i+1}}
 \right)-f\left(\vartheta ,{t_{i}}\right)\right]Y_{t_{i+1}}+\frac{\varepsilon }{\sqrt{\delta
  _\varepsilon} }R_{1,\varepsilon }\\&\qquad+\frac{f\left( {t_{i+1}}
  \right)}{\delta _\varepsilon }\int_{t_{i+1}}^{t_{i+1}+\delta _\varepsilon
}\left[Y_s-Y_{t_{i+1}}\right]\;{\rm d}s\\&\qquad+\frac{1}{\delta _\varepsilon }\int_{t_{i+1}}^{t_{i+1}+\delta _\varepsilon
}\left[f\left( s \right)-f\left( {t_{i+1}}
  \right)\right]Y_s\;{\rm d}s.
\end{align*}
The following estimates are obtained using elementary calculations
\begin{align*}
&\Ex\left(  \left[f\left( {t_{i+1}} \right)-f\left( {t_{i}}
    \right)\right]Y_{t_{i+1}}   \right)^2 \leq C\varphi _\varepsilon ^2,
  \qquad \quad \Ex R_{1,\varepsilon }^2\leq C,\\
&\Ex\left(\frac{f\left( {t_{i+1}}
  \right)}{\delta _\varepsilon }\int_{t_{i+1}}^{t_{i+1}+\delta _\varepsilon
}\left[Y_s-Y_{t_{i+1}}\right]\;{\rm d}s\right)^2\leq C\delta _\varepsilon ,\\
&\Ex\left( \frac{1}{\delta _\varepsilon }\int_{t_{i+1}}^{t_{i+1}+\delta _\varepsilon
}\left[f\left( s \right)-f\left( {t_{i+1}}
  \right)\right]Y_s\;{\rm d}s \right)^2\\
&\qquad \qquad \leq \frac{1}{\delta _\varepsilon
  }\int_{t_{i+1}}^{t_{i+1}+\delta _\varepsilon 
}\left(s-t_{i+1}\right)^2 f'\left( {\tilde s} \right)^2\Ex  Y_s^2\;{\rm d}s\leq C\delta _\varepsilon ^2.
\end{align*}

Therefore
\begin{align*}
&\Ex\left(\frac{X_{t_{i+1}+\delta _\varepsilon }-X_{t_{i+1}}}{\delta _\varepsilon
}-\frac{X_{t_{i}+\delta _\varepsilon }-X_{t_{i}}}{\delta _\varepsilon
}-f\left( {t_{i}}
    \right)\left[Y_{t_{i+1}}-Y_{t_{i}}\right] \right)^2\\
&\qquad\qquad   \leq
C\frac{\varepsilon ^2}{\delta _\varepsilon }+ C\delta _\varepsilon + C\varphi  _\varepsilon^2 .
\end{align*}

These are estimates for the asymptotic derivatives of $X_t$ at the points $t_i
$ and $t_{i+1}$ and  of their difference. 

 The quadratic variation is estimated as follow
\begin{align*}
Q_{\tau ,\varepsilon} &= \sum_{i=0}^{N_{\tau,\varepsilon
}-1} f\left(t_i \right)^2\left(Y_{t_{i+1}}-Y_{t_{i}}
\right)^2=\sum_{i=0}^{N_{\tau,\varepsilon
}-1} f\left(t_i \right)^2\left( \int_{t_{i}}^{t_{i+1}}
a \left( s \right){\rm d}s + \int_{t_{i}}^{t_{i+1}}
b \left( s \right){\rm d}V_s
\right)^2\\
&=\sum_{i=0}^{N_{\tau,\varepsilon
}-1} f\left(t_i \right)^2\left( \int_{t_{i}}^{t_{i+1}}
a \left( s \right){\rm d}s \right)^2+ \sum_{i=0}^{N_{\tau,\varepsilon
}-1} f\left({t_{i}} \right)^2\left( \int_{t_{i}}^{t_{i+1}}
b \left( s \right){\rm d}V_s
\right)^2\\
&\qquad \quad +2\sum_{i=0}^{N_{\tau,\varepsilon
}-1} f\left(t_i \right)^2\int_{t_{i}}^{t_{i+1}}
a \left( s \right){\rm d}s  \int_{t_{i}}^{t_{i+1}}
b \left( s \right){\rm d}V_s,
\end{align*}
where
\begin{align*}
\sum_{i=0}^{N_{\tau,\varepsilon
}-1} f\left(t_i \right)^2\left( \int_{t_{i}}^{t_{i+1}}
a \left( s \right){\rm d}s\right)^2 \leq \varphi _\varepsilon \sum_{i=0}^{N_{\tau,\varepsilon
}-1} f\left(t_i \right)^2 \int_{t_{i}}^{t_{i+1}}
a \left( s \right)^2{\rm d}s \leq C\tau \varphi _\varepsilon  ,
\end{align*}
and
\begin{align*}
\Ex \left(2\sum_{i=0}^{N_{\tau,\varepsilon
}-1} f\left(t_i \right)^2\int_{t_{i}}^{t_{i+1}}
a \left( s \right){\rm d}s  \int_{t_{i}}^{t_{i+1}}
b \left( s \right){\rm d}V_s \right)^2\leq C\tau\varphi _\varepsilon ^2.
\end{align*}
For the stochastic  integral we can write
\begin{align*}
&\sum_{i=0}^{N_{\tau,\varepsilon
}-1} f \left({t_{i}} \right)^2\left(\int_{t_{i}}^{t_{i+1}}
b \left( s \right){\rm d}V_s
\right)^2
\\&\qquad \qquad =\sum_{i=0}^{N_{\tau,\varepsilon
}-1} f  \left({t_{i}}  \right)^2\left(  b \left( {t_{i}}\right)\left[V_{t_{t_{i+1}}}-V_{t_{t_{i}}}\right]+\int_{t_{i}}^{t_{i+1}}
\left[b \left( s \right)-b  \left( {t_{i}}\right)\right]{\rm d}V_s
\right)^2.
\end{align*}
Here
\begin{align*}
&\Ex  \sum_{i=0}^{N_{\tau,\varepsilon
}-1} f\left(t_i \right)^2\left( \int_{t_{i}}^{t_{i+1}}
\left[b \left( s \right)-b\left({t_{i}} \right)\right]{\rm d}V_s
\right)^2\\
&\qquad \qquad =\sum_{i=0}^{N_{\tau,\varepsilon
}-1} f\left(t_i \right)^2 \int_{t_{i}}^{t_{i+1}}
\left[b \left( s \right)-b \left({t_{i}}  \right)\right]^2{\rm d}s\\
&\qquad \qquad =\sum_{i=0}^{N_{\tau,\varepsilon
}-1} f\left({t_{i}} \right)^2 \int_{t_{i}}^{t_{i+1}}\left(s-t_{i}\right)^2
b'\left( {\tilde s} \right)^2{\rm d}s\leq C\tau\varphi _\varepsilon ^2.
\end{align*}
Hence
\begin{align}
\label{23_NE-9}
\Ex\left|Q_{\tau ,\varepsilon }  -\sum_{i=0}^{N_{\tau,\varepsilon
}-1} f \left( t_{i}  \right)^2b \left( {t_{i}} \right)^2
\left[V_{t_{i+1}}-V_{t_{i}}\right]^2\right| \leq C\tau{\varphi _\varepsilon }.
\end{align}
Finally
\begin{align*}
K_{\tau ,\varepsilon} &=\sum_{i=0}^{N_{\tau,\varepsilon
}-1} f \left( t_{i}  \right)^2  b \left( {t_{i}}
\right)^2\frac{\left[V_{t_{t_{i+1}}}-V_{t_{t_{i}}}\right]^2}{t_{i+1}-t_{i}}\left(t_{i+1}-t_{i}\right)\\  
&=\sum_{i=0}^{N_{\tau,\varepsilon
}-1} f \left( t_{i}  \right)^2  b \left( {t_{i}} \right)^2\xi _i^2\left(t_{i+1}-t_{i}\right)\\
&=\sum_{i=0}^{N_{\tau,\varepsilon 
}-1} f \left( t_{i}  \right)^2   b \left( {t_{i}} \right)^2
\left(t_{i+1}-t_{i}\right)+\sum_{i=0}^{N_{\tau,\varepsilon 
}-1} f \left( t_{i}  \right)^2   b \left( {t_{i}} \right)^2\left(\xi
_i^2-1\right)\left(t_{i+1}-t_{i}\right).   
\end{align*}
Here $\xi _i,i=0,\ldots,N_{\tau ,\varepsilon }$ are independent standard
Gaussian random variables, $\xi _i\sim {\cal N}\left(0,1\right)$. Therefore
\begin{align*}
&\Ex\left(K_{\tau ,\varepsilon} -\sum_{i=0}^{N_{\tau,\varepsilon
}-1} f\left({t_{i}}  \right)^2   b \left( {t_{i}} \right)^2 \left(t_{i+1}-t_{i}\right)\right)^2\\
&\qquad \qquad \qquad ={2\varphi _\varepsilon }{}\sum_{i=0}^{N_{\tau,\varepsilon
}-1} f  \left({t_{i}} \right)^4  b  \left({t_{i}}
\right)^4\left(t_{i+1}-t_{i}\right)\leq C\tau\varphi _\varepsilon 
\end{align*}
and
\begin{align}
\label{23_NE-10}
\left|\Psi _\tau -\sum_{i=0}^{N_{\tau,\varepsilon
}-1} f \left({t_{i}} \right)^2  b \left({t_{i}}
\right)^2 \left(t_{i+1}-t_{i}\right)\right|\leq C\tau\varphi _\varepsilon .
\end{align}
The obtained estimates allow us to write
\begin{align*}
\Ex \left|\Psi _{\tau ,\varepsilon} -\sum_{i=0}^{N_{\tau,\varepsilon
}-1} f \left({t_{i}}
\right)^2\left[Y_{t_{i+1}}-Y_{t_{i}}\right]^2\right|\leq C\sqrt{\tau} \left[ \sqrt{\frac{\delta _\varepsilon }{\varphi _\varepsilon }}+\varphi _\varepsilon+\frac{\varepsilon }{\sqrt{\delta _\varepsilon \varphi _\varepsilon }} \right].
\end{align*}
If we put $\delta _\varepsilon =\varepsilon ^q,\varphi _\varepsilon
=\varepsilon ^l$ then the equation
\begin{align*}
\sqrt{\frac{\delta _\varepsilon }{\varphi _\varepsilon }}=\varphi
_\varepsilon=\frac{\varepsilon }{\sqrt{\delta _\varepsilon \varphi
    _\varepsilon }} 
\end{align*}
gives us the values $q=1$ and $l=1/3$, i.e., $\delta _\varepsilon =\varepsilon
$ and $\varphi _\varepsilon =\varepsilon ^{1/3}$. Hence
\begin{align*}
\Ex\left|\Psi _{\tau ,\varepsilon }-\sum_{i=0}^{N_{\tau,\varepsilon
}-1} f \left({t_{i}}
\right)^2\left[Y_{t_{i+1}}-Y_{t_{i}}\right]^2\right|\leq C\,{\tau }^{1/2}\,\varepsilon
^{1/3}. 
\end{align*}
This estimate together with \eqref{23_NE-9} and \eqref{23_NE-10} allows us to write
\begin{align*}
\Ex_\vartheta \left|\Psi _{\tau ,\varepsilon} -\Psi _\tau \right|\leq
C\,\tau^{1/2} \,\varepsilon ^{1/3}.
\end{align*}

The similar calculations provide the estimate of the  polynomial moments too.

\end{proof}

Remark that it is possible to prove the asymptotic normality
\begin{align*}
\tau ^{-1/2}\varepsilon ^{-1/3}\left(\Psi _{\tau ,\varepsilon} -\Psi _\tau
\right)\Longrightarrow {\cal N}\left(0,D^2\right)
\end{align*}
with some limit variance $D^2 $. The estimate \eqref{23_NE-8} is an intermediate
result, which will be used in the next section for the construction of the
preliminary estimator for One-step MLE-process.

\begin{remark}
\label{R23_NE-1}
{\rm The obtained rate $\varepsilon ^{1/3}$ is not  optimal but is
  sufficient for the construction of the preliminary estimator below.   An interesting result concerning estimation of the function $\Psi _\tau $
  by observations  \eqref{23_NE-6}, \eqref{23_NE-7} was obtained in \cite{R11}.
    Supposing that 
  $a\left(t\right)\equiv 0$, $\sigma \left(t\right)\equiv \sigma $, $f
  \left(t\right)\equiv 1$ it was shown that this
  model of observations 
is asymptotically equivalent (in Le Cam's sense) to the model of observations
\begin{align*}
{\rm d}X_t=\sqrt{2b\left(t\right)}\,{\rm d}t+\sigma ^{1/2}\varepsilon ^{1/2}{\rm d}W_t.
\end{align*}
This equivalence  shows why the rate $\varepsilon ^{1/2} $ is  ``more natural'' for
the model \eqref{23_NE-6}, \eqref{23_NE-7}. 

We suppose that more detailed analysis of estimator $\Psi _{\tau ,\varepsilon
}$ can show the rate $\varepsilon ^{1/2} $, but for further application of
$\Psi _{\tau ,\varepsilon }$ in the construction of One-step MLE process  the rate $\varepsilon ^{1/3} $ is sufficient.
}
\end{remark}

\section{Method of moments estimator}

Let us return to the parametric   partially observed linear system
\begin{align}
\label{23_MME-11}
{\rm d}X_t&=f\left(\vartheta ,t\right)Y_t\,{\rm d}t+\varepsilon \sigma
\left(t\right){\rm d}W_t,\qquad  X_0=0,\qquad 0\leq t\leq T,\\
{\rm d}Y_t&=a\left(\vartheta,t\right)Y_t\,{\rm d}t+b\left(\vartheta ,t\right){\rm
  d}V_t,\qquad \;Y_0=y_0
\label{23_MME-12}
\end{align}
The process $X^T=\left(X_t,0\leq t\leq T\right)$ is observed, the process
$Y^T=\left(Y_t,0\leq t\leq T\right)$ is {\it hidden} and the Wiener processes
$V^T=\left(V_t,0\leq t\leq T\right)$ and $W^T=\left(W_t,0\leq t\leq T\right)$
as before are independent.  The functions $b\left(\cdot \right)$ and
$f\left(\cdot \right)$ are supposed to be known. The functions $a\left(\cdot
\right)$ and $\sigma \left(\cdot \right)$ do not used in the construction of
the estimator. The parameter $\vartheta\in \Theta
=\left(\alpha ,\beta \right) $ is unknown and has to be estimated by
observations $X^T$.

 Fix some $\tau\in (0,T] $ and define the function
\begin{align*}
\Psi _t \left(\vartheta \right)=\int_{0}^{t }f\left(\vartheta
,s\right)^2b\left(\vartheta ,s\right)^2{\rm d}s ,\qquad 0<t\leq \tau ,\qquad
\vartheta \in \Theta\subset{\cal R}^d .
\end{align*}
As estimator of this function we take the nonparametric estimator 
\begin{align*}
\Psi _{t ,\varepsilon} =\sum_{i=0}^{N_{t,\varepsilon
}-1}\left(\frac{X_{t_{i+1}+\varepsilon }-X_{t_{i+1}}}{\varepsilon }-\frac{X_{t_{i}+\varepsilon }-X_{t_{i}}}{\varepsilon }\right)^2,\qquad\quad  0<t\leq \tau.  
\end{align*}
Here $t_{i+1}-t_i=\varepsilon ^{1/3}, t_0=0$ and $N_{t,\varepsilon
}=\left[t\varepsilon ^{-1/3}\right]$.

The MME $\check\vartheta _{\tau ,\varepsilon }$ can be defined as follows. Let
us fix $d$ points $0<\tau _1<\tau _2<\ldots<\tau _d=\tau$ and denote the
vectors $\Psi_{\tau ,\varepsilon}^*=\left(\Psi_{\tau_1
  ,\varepsilon},\Psi_{\tau_2 ,\varepsilon},\ldots,\Psi_{\tau_d
  ,\varepsilon}\right)^\top $ and $\Psi_{\tau }^*\left(\vartheta
\right)=\left(\Psi_{\tau_1 }\left(\vartheta \right),\Psi_{\tau_2
}\left(\vartheta \right),\right.$ $ \left.\ldots,\Psi_{\tau_d }\left(\vartheta
\right)\right)^\top $. The estimator is defined as solution of the system of
equations
\begin{align}
\label{23_MME-13}
\Psi_{\tau ,\varepsilon}^*=\Psi_{\tau }^*(\check\vartheta _{\tau
  ,\varepsilon } ).
\end{align}
If this equation has no solution then  we put $\check\vartheta _{\tau
  ,\varepsilon }=\vartheta ^\circ $, where $\vartheta ^\circ $ is some value
$\not\in \bar\Theta $. 

Note that having the estimator  $\Psi _{t ,\varepsilon},0\leq t\leq \tau  $ we
can define the MDE $\check\vartheta _{\tau ,\varepsilon }$ too  by the equation
\begin{align}
\label{23_MME-14}
\int_{0}^{\tau }\left[\Psi _{t ,\varepsilon}-\Psi _t (\check\vartheta
  _{\tau ,\varepsilon } )\right]^2{\rm d}t=\inf_{\vartheta \in\Theta
}\int_{0}^{\tau }\left[\Psi _{t ,\varepsilon}-\Psi _t \left(\vartheta
  \right)\right]^2{\rm d}t.
\end{align}
If these equations has more than one solution then anyone of them can be taken
as the estimator of $\vartheta $.

Introduce the notation
\begin{align*}
g\left(\vartheta _0,\nu \right)=\inf_{\left\|\vartheta -\vartheta
  _0\right\|\geq \nu }\left\|\Psi _\tau  \left(\vartheta 
  \right)-\Psi _\tau  \left(\vartheta_0
  \right)\right\|,\qquad {\bf T}\left(\vartheta _0\right)=\dot \Psi_{\tau }^*(
\vartheta _0 ) \dot \Psi_{\tau }^*(\vartheta _0 )^\top
\end{align*}
 
\begin{proposition}
\label{P23_MME-2} Suppose that the following  conditions hold
\begin{itemize}

\item  The functions $a\left(\vartheta ,t
  \right),b\left(\vartheta ,t \right),f\left(\vartheta ,t \right), \sigma \left(t
  \right), \left(\vartheta ,t\right)\in  \bar\Theta \times\left[0,T\right] $
  have continuous deri- vatives w.r.t. $t$ and two continuous derivatives
  w.r.t. $\vartheta $. 

\item  The matrix ${\bf T}\left(\vartheta _0\right) $ is uniformly non
  degenerate
\begin{align*}
\inf_{\vartheta _0\in\Theta }\;\;\inf_{\left\|e\right\|=1,e\in {\cal R}^d} e^\top
    {\bf T}\left(\vartheta _0\right)e>0. 
\end{align*}
\item For any $\nu >0$
\begin{align*}
\inf_{\vartheta _0\in\Theta}g\left(\vartheta _0,\nu \right)>0.
\end{align*}
\end{itemize}
 Then the MME $\check
\vartheta _{\tau ,\varepsilon } $ is uniformly   consistent
and for any $p\geq 2$ there exists a constant $C=C\left(p\right)>0$ such that
\begin{align}
\label{23_MME-15}
\sup_{\vartheta _0\in\Theta }\varepsilon ^{-p/3}\,\tau ^{p/2}\Ex_{\vartheta _0}\left\|\check
\vartheta _{\tau ,\varepsilon }-\vartheta _0\right\|^p\leq C .
\end{align}
\end{proposition}
\begin{proof}
Remark that if $\check\vartheta _{\tau ,\varepsilon }\in\Theta $ then it
satisfies the equation
\begin{align*}
\inf_{\vartheta \in\Theta }\left\|\Psi_{\tau ,\varepsilon}^*-\Psi_{\tau }^*(\vartheta )\right\|=\left\|\Psi_{\tau ,\varepsilon}^*-\Psi_{\tau }^*(\check\vartheta _{\tau
  ,\varepsilon } )\right\|=0.
\end{align*}
Below we use elementary inequalities and the estimate \eqref{23_NE-8}
\begin{align*}
&\sup_{\vartheta _0\in\Theta }\Pb_{\vartheta _0}\left(\left\|\check\vartheta _{\tau ,\varepsilon
}-\vartheta_0 \right\|>\nu \right)\\
&\qquad \qquad =\sup_{\vartheta _0\in\Theta }\Pb_{\vartheta
  _0}\left(\inf_{\left\|\vartheta-\vartheta _0\right\|<\nu }\left\|\Psi_{\tau
  ,\varepsilon}^*-\Psi_{\tau }^*(\vartheta )\right\|>
\inf_{\left\|\vartheta-\vartheta _0\right\|\geq \nu }\left\|\Psi_{\tau
  ,\varepsilon}^*-\Psi_{\tau }^*(\vartheta )\right\|\right)\\ 
&\qquad \qquad \leq
\sup_{\vartheta _0\in\Theta }\Pb_{\vartheta _0}\left(\left\|\Psi_{\tau ,\varepsilon}^*-\Psi_{\tau
}^*(\vartheta_0 )\right\|+\inf_{\left\|\vartheta-\vartheta _0\right\|<\nu
}\left\|\Psi_{\tau }^*\left(\vartheta \right)-\Psi_{\tau }^*(\vartheta_0
)\right\|\right.\\
 &\qquad \qquad\qquad \left. > \inf_{\left\|\vartheta-\vartheta
  _0\right\|\geq \nu }\left\|\Psi_{\tau }^*\left(\vartheta \right)-\Psi_{\tau
}^*(\vartheta_0 )\right\|-\left\|\Psi_{\tau ,\varepsilon}^*-\Psi_{\tau
}^*(\vartheta_0 )\right\|\right)\\ 
&\qquad \qquad=\sup_{\vartheta _0\in\Theta }\Pb_{\vartheta
  _0}\left(2\left\|\Psi_{\tau ,\varepsilon}^*-\Psi_{\tau }^*(\vartheta_0
)\right\| > g\left(\vartheta _0,\nu \right)\right)\\
 &\qquad \qquad\leq
\frac{2^{2N}}{\inf_{\vartheta _0\in\Theta }g\left(\vartheta _0,\nu \right)^{2N}}\sup_{\vartheta _0\in\Theta }\Ex_{\vartheta
  _0}\left\|\Psi_{\tau ,\varepsilon}^*-\Psi_{\tau }^*(\vartheta_0
)\right\|^{2N}\\
 &\qquad \qquad\leq
\frac{C\tau ^N\varepsilon ^{2N/3}}{\inf_{\vartheta _0\in\Theta }g\left(\vartheta _0,\nu \right)^{2N}}\rightarrow 0.
\end{align*}
Therefore the estimator $\check\vartheta _{\tau ,\varepsilon
} $ is uniformly  consistent.

Recall that the MME  $\check\vartheta _{\tau ,\varepsilon
}$ is one of the solutions of the MMEq (method of moments equation)
\begin{align*}
\left(\Psi_{\tau ,\varepsilon}^*-\Psi_{\tau }^*(\check
\vartheta _{\tau ,\varepsilon } )\right)^\top \dot \Psi_{\tau }^*(\check
\vartheta _{\tau ,\varepsilon } )=0.
\end{align*}
Using the consistency of this estimator we can write 
\begin{align*}
0&=\left(\Psi_{\tau ,\varepsilon}^*-\Psi_{\tau }^*(
\vartheta _0 )+\Psi_{\tau }^*(
\vartheta _0 )  - \Psi_{\tau }^*(\check
\vartheta _{\tau ,\varepsilon } )\right)^\top \dot \Psi_{\tau }^*(\check
\vartheta _{\tau ,\varepsilon } )\\
&=\left(\Psi_{\tau ,\varepsilon}^*-\Psi_{\tau }^*(
\vartheta _0 )-  \dot \Psi_{\tau }^*(\check
\vartheta _0 )^\top (\check
\vartheta _{\tau ,\varepsilon }-\vartheta _0 )+ O\left(\varepsilon
^{2/3}\right)\right)^\top \left(\dot \Psi_{\tau }^*(
\vartheta _0 )+ O\left(\varepsilon ^{1/3}\right)\right)\\ 
&=\left(\Psi_{\tau ,\varepsilon}^*-\Psi_{\tau }^*(
\vartheta _0 )\right)^\top \dot \Psi_{\tau }^*(
\vartheta _0 )     -  (\check
\vartheta _{\tau ,\varepsilon }-\vartheta _0 )\dot \Psi_{\tau }^*(\check
\vartheta _0 )^\top \dot \Psi_{\tau }^*(
\vartheta _0 )   + O\left(\varepsilon
^{2/3}\right) .
\end{align*}
Therefore
\begin{align*}
\varepsilon ^{-1/3}(\check \vartheta _{\tau ,\varepsilon }-\vartheta _0 )&=
\left(\dot \Psi_{\tau }^*(\check \vartheta _0 )^\top \dot \Psi_{\tau }^*( \vartheta
_0 ) \right)^{-1}\;\varepsilon ^{-1/3} \left(\Psi_{\tau
  ,\varepsilon}^*-\Psi_{\tau }^*( \vartheta _0 )\right)^\top \dot \Psi_{\tau
}^*( \vartheta _0 ) +O\left(\varepsilon ^{1/3}\right)\\
 &={\bf  T}\left(\vartheta _0\right)^{-1} \;\varepsilon ^{-1/3} \left(\Psi_{\tau
  ,\varepsilon}^*-\Psi_{\tau }^*( \vartheta _0 )\right)^\top \dot \Psi_{\tau
}^*( \vartheta _0 ) +O\left(\varepsilon ^{1/3}\right) .
\end{align*}
As the matrix ${\bf  T}\left(\vartheta _0\right) $ is uniformly non degenerate
and we have \eqref{23_NE-8} these allow us to write the estimate: for any $p\geq 2$
\begin{align*}
\varepsilon ^{-p/3}\sup_{\vartheta _0\in\Theta }\Ex_{\vartheta
  _0}\left\|\check \vartheta _{\tau ,\varepsilon }-\vartheta _0\right\|^{p}
\leq C\,\varepsilon ^{-p/3}\sup_{\vartheta _0\in\Theta }\Ex_{\vartheta
  _0}\left\|\Psi_{\tau ,\varepsilon}^*-\Psi_{\tau }^*( \vartheta _0 )
\right\|^{p}\leq C\tau ^{p/2}
\end{align*}

\end{proof}

\begin{example}
\label{E23-1} {\rm  Suppose that we have the model of observations
\eqref{23_MME-11},\eqref{23_MME-12}, where $f\left(\vartheta
,t\right)=\vartheta f\left(t\right), \vartheta \in\left(\alpha ,\beta \right),
\alpha >0$, $b\left(\vartheta ,t\right)=b\left(t\right)$ and all corresponding
conditions are fulfilled. Then
\begin{align*}
\Psi _\tau \left(\vartheta \right)=\vartheta ^2\int_{0}^{\tau
}f\left(t\right)^2b\left(t\right)^2{\rm d}t 
\end{align*}
and the MME is
\begin{align*}
\check\vartheta _{\tau ,\varepsilon }=\sqrt{\Psi _{\tau,\varepsilon  }}\left(\int_{0}^{\tau
}f\left(t\right)^2b\left(t\right)^2{\rm d}t \right)^{-1/2}.
\end{align*}
This estimator is consistent and has rate of convergence ${\varepsilon }^{1/3}$
(Proposition \ref{P23_MME-2}). 
}
\end{example}
\bigskip

\begin{example}
\label{E23-2} {\rm
Consider the model \eqref{23_MME-11},\eqref{23_MME-12}, where
$f\left(\vartheta, t\right)=f\left(t\right)$ and
\begin{align*}
b\left(\vartheta
,t\right)=\left({h\left(t\right)+\sum_{k=1}^{d}\vartheta_k
  g_k\left(t\right)}\right)^{1/2} .
\end{align*}
 Here the
functions $h\left(\cdot \right)$ and $g_k\left(\cdot \right),k=1,\ldots,d$ are
positive and all $\vartheta  >0$. Then
\begin{align*}
\Psi _{\tau _k} \left(\vartheta \right)&=\int_{0}^{\tau _k
}f\left(s\right)^2h\left(s\right){\rm d}s +\sum_{j=1}^{d}\vartheta_j
\int_{0}^{\tau _k }f\left(s\right)^2 g_j\left(s\right) {\rm d}s\\
&=H_k+\sum_{j=1}^{d}\vartheta_jG_{jk},\qquad \quad k=1,\ldots,d,
\end{align*}
and we have the relation
\begin{align*}
\vartheta ={\bf G}^{-1}   \bigl(\Psi _{\tau } \left(\vartheta \right)-H\bigr)
\end{align*}
with the obvious notations. 

Hence the MME is
\begin{align*}
\check\vartheta _{\tau ,\varepsilon }={\bf G}^{-1} \bigl(\Psi _{\tau,\varepsilon }-H\bigr)
\end{align*}
and this estimator has the properties described in the Proposition \ref{P23_MME-2}.
}
\end{example}

\section{One-step MLE-process.}

 We have the same  partially observed linear system
\begin{align}
\label{23-56}
{\rm d}X_t&=f\left(\vartheta ,t\right)Y_t{\rm d}t+\varepsilon \sigma
\left(t\right){\rm d}W_t,\qquad \quad   X_0=0,\qquad 0\leq t\leq T,\\
{\rm d}Y_t&=a\left(\vartheta ,t\right)Y_t{\rm d}t+b\left(\vartheta ,t\right){\rm
  d}V_t,\qquad \quad Y_0=0.\nonumber
\end{align}
As before the process $X^T=\left(X_t,0\leq t\leq T\right)$ is observable and
$Y^T$ is hidden.  All functions $a\left(\cdot \right),b\left(\cdot
\right),f\left(\cdot \right)$ and $\sigma \left(\cdot \right)$ are supposed to
be known. The parameter $\vartheta\in \Theta \subset{\cal R}^d $ is
unknown and has to be estimated by observations $X^T$.

Remind that the MLE $\hat \vartheta _\varepsilon $ and BE $\tilde\vartheta
_\varepsilon $ are  defined by the relations \eqref{22-45a} and their
calculation requires solutions  $m\left(\vartheta ,\cdot \right),\vartheta
\in\Theta  $, $\gamma \left(\vartheta ,\cdot \right),\vartheta
\in\Theta $ of the equations 
\begin{align}
\label{23-57}
{\rm d}m\left(\vartheta ,t\right)&=a\left(\vartheta ,t\right)m\left(\vartheta
,t\right){\rm d}t +\frac{\gamma \left(\vartheta ,t\right)f\left(\vartheta
  ,t\right)}{\varepsilon ^2\sigma \left(t\right)^2}\left[{\rm
    d}X_t-f\left(\vartheta ,t\right)m\left(\vartheta ,t\right){\rm d}t\right],\;m\left(\vartheta ,0\right)=y_0
\end{align}
and
\begin{align}
\label{23-58}
\frac{\partial \gamma \left(\vartheta ,t\right)}{\partial
  t}&=2a\left(\vartheta ,t\right)\gamma \left(\vartheta ,t\right)- \frac{\gamma
  \left(\vartheta ,t\right)^2f \left(\vartheta ,t\right)^2}{\varepsilon
  ^2\sigma \left(t\right)^2} +b\left(\vartheta ,t\right)^2 ,\qquad \gamma
\left(\vartheta ,0\right)=0.  
\end{align}

 Our goal is to construct an approximation $\hat m_t$ of the process
 $m\left(\vartheta ,t\right)$ (adaptive filter) by replacing $\vartheta $ in
 the equations \eqref{23-57}, \eqref{23-58} by some estimator-process
 $\bar\vartheta _{t,\varepsilon }, 0<t\leq T$. Therefore we need a sequence of
 estimators $\bar\vartheta _{t,\varepsilon },0<t\leq T$ which are measurable
 w.r.t. observations up to the time $t$ for each $t\in (0,T]$.  To avoid the
   evident numerical difficulties related with the calculations MLE
   $\hat\vartheta _{t,\varepsilon },0<t\leq T$ we consider an approach based
   on the Fisher-score device and MME $\check \vartheta _{\tau ,\varepsilon }$
   as preliminary one.   The proposed
   One-step MLE-process allows us to obtain an estimator-process, which is easy
   to calculate, has good asymptotic properties and can be used for the
   construction of adaptive filter.

 The Fiser-score function for this model is
\begin{align*}
\frac{\partial }{\partial \vartheta }\ln L\left(\vartheta
,X^T\right)=\int_{0}^{T}\frac{\dot M\left(\vartheta ,t\right)}{\varepsilon ^2\sigma \left(t\right)^2}\left[{\rm d}X_t-M\left(\vartheta ,t\right){\rm d}t\right],
\end{align*}
where $M\left(\vartheta ,t\right)=f\left(\vartheta ,t\right)\left(\vartheta
,t\right)\left(\vartheta ,t\right)m$. 
Let us consider the construction of One-step MLE-process in the case of
observations \eqref{23-56}.  Suppose that for some fixed small value $\tau >0$
we already have MME $\check \vartheta _{\tau ,\varepsilon }$ constructed by
the observations $X^\tau =\left(X_t,0\leq t\leq \tau \right)$. The solution of
the equation \eqref{23-57} on the time interval $\left[0,\tau \right]$ can be
written as follows
\begin{align*}
m\left(\vartheta ,t\right)&=y_0\Phi\left(\vartheta,0,
t\right)+ \Phi\left(\vartheta,0,
t\right)\int_{0}^{t}\Phi\left(\vartheta,0,
s\right)^{-1} B\left(\vartheta ,s\right){\rm d}X_s\\
& = y_0\Phi\left(\vartheta,0,
t\right)+ \Phi\left(\vartheta,0,
t\right)\int_{0}^{t}H\left(\vartheta ,s\right){\rm d}X_s,
\end{align*}
where 
\begin{align*}
\Phi\left(\vartheta,0,
t\right)&=\exp\left(\int_{0}^{t}A
\left(\vartheta ,v\right){\rm d}v \right) ,\qquad \quad A
\left(\vartheta ,t\right)=a\left(\vartheta ,t\right)- \frac{\gamma
  \left(\vartheta ,t\right)f\left(\vartheta ,t\right)^2}{\varepsilon ^2\sigma
  \left(t\right)^2},\\ 
B \left(\vartheta ,s\right)&= \frac{\gamma
  \left(\vartheta ,s\right)f\left(\vartheta ,s\right)}{\varepsilon ^2\sigma
  \left(s\right)^2},\qquad \quad \qquad H\left(\vartheta
,s\right)=\Phi\left(\vartheta,0, s\right)^{-1} \frac{\gamma \left(\vartheta
  ,s\right)f\left(\vartheta ,s\right)}{\varepsilon ^2\sigma \left(s\right)^2}.
\end{align*}
We can not put $\check\vartheta _{\tau,\varepsilon  } $ in $m\left(\vartheta
,t\right) $ because the corresponding stochastic integral It\^o is not
defined. Note that the function $H\left(\vartheta ,s\right)$ is continuously
differentiable w.r.t. $s$, hence the following relation
\begin{align*}
\int_{0}^{\tau }H\left(\vartheta ,s\right){\rm d}X_s=X_\tau H\left(\vartheta
,\tau \right)-\int_{0}^{\tau }X_sH'\left(\vartheta ,s\right){\rm d}s
\end{align*}
holds. Now we can put
\begin{align}
\label{23-59}
m\left(\check\vartheta _{\tau,\varepsilon  } ,\tau
\right)&=y_0\Phi\left(\check\vartheta _{\tau,\varepsilon  },0,
\tau \right)+\Phi\left(\check\vartheta _{\tau,\varepsilon  } ,0,\tau \right)
X_\tau H\left(\check\vartheta _{\tau,\varepsilon  } 
,\tau \right)\nonumber\\
&\qquad -\Phi\left(\check\vartheta _{\tau,\varepsilon  } ,0,\tau \right)\int_{0}^{\tau }X_sH'\left(\check\vartheta _{\tau,\varepsilon  } ,s\right){\rm d}s .
\end{align}
Let us define the random process $\hat m\left(t\right),\tau \leq t\leq T$ by
the equation
\begin{align}
\label{23-60}
{\rm d}\hat m\left(t\right)&=a\left(\check\vartheta _{\tau,\varepsilon }
,t\right)\hat m\left(t\right){\rm d}t +\frac{\gamma
  \left(\check\vartheta _{\tau,\varepsilon } ,t\right)f\left(\check\vartheta
  _{\tau,\varepsilon } ,t\right)}{\varepsilon ^2\sigma
  \left(t\right)^2}\left[{\rm d}X_t-f\left(\check\vartheta _{\tau,\varepsilon
  } ,t\right)\hat m\left(t\right){\rm d}t\right]
\end{align}
with the initial value $\hat m\left(\tau \right)=m\left(\check\vartheta
_{\tau,\varepsilon } ,\tau \right)$. The solution of this equation is
\begin{align*}
\hat m\left(t\right)=m\left(\check\vartheta
_{\tau,\varepsilon } ,\tau \right)\Phi\left(\check\vartheta
  _{\tau,\varepsilon },\tau ,t\right)+\int_{\tau }^{t}\Phi\left(\check\vartheta
  _{\tau,\varepsilon },s ,t\right)\frac{\gamma
  \left(\check\vartheta _{\tau,\varepsilon } ,s\right)f\left(\check\vartheta
  _{\tau,\varepsilon } ,s\right)}{\varepsilon ^2\sigma
  \left(s\right)^2} {\rm d}X_s.
\end{align*}

We need as well the equation for the derivative $\dot m\left(\vartheta
,t\right)$ on the interval $\left[\tau ,T\right]$ at the point $\vartheta
=\check\vartheta _{\tau,\varepsilon }$ (denoted $\dot m\left(\check\vartheta
_{\tau,\varepsilon } ,t\right)= \hat{\dot m}\left(t\right) $), which we obtain
by formal differentiation of the equation \eqref{23-57}
\begin{align}
\label{23-61}
{\rm d}\hat{\dot m}\left(t\right)=A\left(\check\vartheta _{\tau,\varepsilon }
,t\right)\hat{\dot m}\left(t\right){\rm d}t+\dot A
\left(\check\vartheta _{\tau,\varepsilon } ,t\right) \hat{ m}\left(t\right){\rm d}t+\dot B\left(\check\vartheta
_{\tau,\varepsilon } ,t\right){\rm d}X_t,\quad \tau \leq t\leq T
\end{align}
with the initial value $\hat{\dot m}\left(\tau \right)=\dot
m\left(\check\vartheta _{\tau,\varepsilon },\tau \right) $, which can be
calculated similar to the given above representation of the corresponding
stochastic integral like \eqref{23-59}.

 Introduce notations:  
\begin{align*}
M\left(\check\vartheta _{\tau,\varepsilon }
 ,t\right)&=f\left(\check\vartheta _{\tau,\varepsilon },t \right)\hat
 m\left(t\right),\qquad \quad \dot M\left(\check\vartheta _{\tau,\varepsilon }
 ,t\right)=\dot f\left(\check\vartheta _{\tau,\varepsilon },t \right)\hat
 m\left(t\right)+f\left(\check\vartheta _{\tau,\varepsilon },t \right)\hat{\dot m}
 \left(t\right),\\
{\bf I}_\tau ^t\left(\vartheta \right)&=
\int_{\tau }^{t}\frac{\dot S\left(\vartheta ,s\right)\dot S\left(\vartheta ,s\right)^\top }{2S\left(\vartheta
  ,s\right)\sigma \left(s\right)}{\rm d}s,\qquad\qquad  \eta   _{t,\varepsilon }=\frac{\vartheta _{t,\varepsilon
}^\star-\vartheta _0}{\sqrt{\varepsilon }} ,
\\
\vartheta _{t,\varepsilon }^\star&=\check\vartheta _{\tau,\varepsilon
}+{{\bf I}_\tau^t (\check\vartheta _{\tau,\varepsilon  } )}^{-1} \int_{\tau
}^{t}\frac{\dot M(\check\vartheta _{\tau,\varepsilon  },s)}{{\varepsilon}
  \sigma
    \left(s\right)^2} \left[{\rm d}X_s-M(\check\vartheta _{\tau,\varepsilon
  },s ){\rm d}s\right],\qquad \tau <t\leq T,\\
  \eta _t\left(\vartheta _0\right)&={{\bf I}_\tau^t
  (\vartheta _{0  } )}^{-1}\int_{\tau 
}^{t}\frac{\dot S\left(\vartheta _0,s\right)}{\sqrt{2S\left(\vartheta
  _0,s\right)\sigma \left(s\right)}}\;{\rm d}w\left(s \right)  ,\qquad \quad \tau
<t\leq T,
\end{align*}
where $w\left(s\right),\tau \leq s\leq T$ is some Wiener process.

\begin{theorem}
\label{T23-1} Let the following  assumptions hold.
\begin{enumerate}

\item The functions $f\left(\vartheta ,t \right),b\left(\vartheta ,t
  \right),\sigma \left(t \right),\left(\vartheta ,t\right)\in \Theta
  \times\left[0,T\right] $ are separated from zero.

\item The conditions of Proposition \ref{P23_MME-2} are fulfilled.

\item For any $t_0\in (\tau ,T]$
\begin{align*}
\inf_{\vartheta_0 \in\Theta } \quad  \inf_{\left\|e\right\|=1,e\in{\cal R}^d }
e^\top  {\bf I}_\tau ^{t_0}\left(\vartheta_0 \right)e >0.
\end{align*}

\end{enumerate}
Then the One-step MLE-process $\vartheta _{t,\varepsilon }^\star,\tau <t\leq T
$ has the following properties:
\begin{enumerate}
\item It  is uniformly consistent: for any $t_0\in (\tau ,T] $  and any   $\nu >0$
\begin{align*}
\Pb_{\vartheta _0}\left(\sup_{t_0\leq t\leq T} \left\|\vartheta _{t,\varepsilon
}^\star-\vartheta _0 \right\|\geq \nu \right) \longrightarrow 0,
\end{align*}
\item The random process $\eta _{t,\varepsilon },t_0\leq t\leq T$ converges in
  distribution in the measurable space $\left({\cal C}\left[t_0,T],{\scr
    B}\right]\right)$ to the Gaussian process $\eta _{t }\left(\vartheta
    _0\right),t_0\leq t\leq T$
\begin{align*}
\eta _{\cdot ,\varepsilon }\Longrightarrow \eta _{\cdot }\left(\vartheta
_0\right),\qquad \qquad \eta _{t }\left(\vartheta _0\right)\sim {\cal
  N}\left(0,{\bf I}_\tau ^{t}\left(\vartheta_0 \right)^{-1}\right),
\end{align*}
\item All polynomial moments converge: for any $p>0$

\begin{align}
\label{23-62}
\varepsilon ^{-p/2}\Ex_{\vartheta _0}\|\vartheta _{t,\varepsilon
}^\star-\vartheta _0\|^p\longrightarrow \Ex_{\vartheta _0}\|\eta _{t
}\left(\vartheta_0 \right)\|^p. 
\end{align}
\end{enumerate}
\end{theorem}
\begin{proof}
The uniform consistency follows from the representation 
\begin{align*}
\sup_{t_0\leq t\leq T}\left\|\vartheta _{t,\varepsilon} ^\star-\vartheta
_0\right\|&\leq \left\|\check\vartheta _{\tau,\varepsilon }-\vartheta
_0\right\|+C\varepsilon \sup_{t_0\leq t\leq T}\left\|\int_{\tau _\varepsilon
}^{t}\frac{{\dot M} \left(\check\vartheta _{\tau,\varepsilon
  },s\right)}{\sigma \left(s\right)}{\rm d}W_s\right\|\\
 &\qquad +C\int_{\tau  _\varepsilon }^{T}{\left\|{\dot M} \left(\check\vartheta _{\tau,\varepsilon
  },s\right)\right\|}\left|M \left(\vartheta _0,s\right) -M
\left(\check\vartheta _{\tau ,\varepsilon },s\right)\right|{\rm d}s 
\end{align*}
and consistency of the MME $\check\vartheta _{\tau,\varepsilon }$.

We have to study the following expression 
\begin{align*}
\frac{\vartheta _{t,\varepsilon }^\star-\vartheta _0}{\sqrt{\varepsilon
}}&=\frac{\left(\check\vartheta _{\tau,\varepsilon }-\vartheta _0
  \right)}{\sqrt{\varepsilon }} 
+{{\bf I}_\tau^t (\check\vartheta _{\tau,\varepsilon } )}^{-1} \int_{\tau
}^{t}\frac{\dot M(\check\vartheta _{\tau,\varepsilon },s)}{\sqrt{\varepsilon}
  \sigma \left(s\right)} {\rm d}\bar W_s\\
&\qquad  +{{\bf I}_\tau^t (\check\vartheta _{\tau,\varepsilon } )}^{-1}\int_{\tau
}^{t}\frac{\dot M(\check\vartheta _{\tau,\varepsilon },s)}{\varepsilon^{3/2}
  \sigma \left(s\right)^2} \left[f\left(\vartheta _0,s\right)m\left(\vartheta 
  _0,s\right)-M(\check\vartheta 
  _{\tau,\varepsilon },s )\right]{\rm d}s.
\end{align*}
We have the following relations
\begin{align*}
&{{\bf I}_\tau^t (\check\vartheta _{\tau,\varepsilon } )}^{-1}={{\bf I}_\tau^t
  (\vartheta _0)}^{-1}+O\left(\left\|\check\vartheta _{\tau,\varepsilon
}-\vartheta _0\right\|\right),\\
&f\left(\vartheta _0,s\right)m\left(\vartheta _0,s\right)-M(\check\vartheta
  _{\tau,\varepsilon },s )=\left(\vartheta _0-\check\vartheta
  _{\tau,\varepsilon }\right)^\top \dot M\left(\vartheta
  _0,t\right)+O\left(\left\|\check\vartheta _{\tau,\varepsilon 
}-\vartheta _0\right\|\right),\\
& \int_{\tau
}^{t}\frac{\dot M(\check\vartheta _{\tau,\varepsilon },s)}{\sqrt{\varepsilon}
  \sigma \left(s\right)} {\rm d}\bar W_s=\int_{\tau
}^{t}\frac{\dot M(\vartheta _0,s)}{\sqrt{\varepsilon}
  \sigma \left(s\right)} {\rm d}\bar W_s+O\left(\left\|\check\vartheta
  _{\tau,\varepsilon }-\vartheta _0\right\|\right).
\end{align*}
Therefore
\begin{align*}
&{{\bf I}_\tau^t (\check\vartheta _{\tau,\varepsilon } )}^{-1}\int_{\tau
  }^{t}\frac{\dot M(\check\vartheta _{\tau,\varepsilon },s)}{\varepsilon^{3/2}
    \sigma \left(s\right)^2} \left[f\left(\vartheta
    _0,s\right)m\left(\vartheta _0,s\right)-M(\check\vartheta
    _{\tau,\varepsilon },s )\right]{\rm d}s\\
&\qquad \qquad ={{\bf I}_\tau^t (\vartheta _0)}^{-1}\int_{\tau
}^{t}\frac{\dot M(\vartheta _0,s)\dot M(\vartheta _0,s)^\top }{\varepsilon^{3/2}
  \sigma \left(s\right)^2}{\rm d}s\left(\vartheta _0-\check\vartheta
  _{\tau,\varepsilon }\right)+O\left(\left\|\check\vartheta _{\tau,\varepsilon
}-\vartheta _0\right\|\right).
\end{align*}

The random process $M\left(\vartheta ,t\right)$ is solution of the  stochastic
differential equation
\begin{align*}
{\rm d}M\left(\vartheta ,t\right)&=f'\left({\vartheta ,t}\right)m\left(\vartheta
,t\right){\rm d}t+a\left(\vartheta ,t\right)M\left(\vartheta ,t\right){\rm
  d}t-\frac{\gamma \left(\vartheta ,t\right)f\left(\vartheta
  ,t\right)^2}{\varepsilon ^2\sigma \left(t\right)^2}M\left(\vartheta
,t\right){\rm d}t\\
&\qquad \qquad \qquad  +\frac{\gamma \left(\vartheta ,t\right)f\left(\vartheta
  ,t\right)^2}{\varepsilon ^2\sigma \left(t\right)^2}{\rm d}X_t\\
&=-q_\varepsilon \left(\vartheta ,t\right)M\left(\vartheta ,t\right){\rm d}t+\frac{\gamma \left(\vartheta ,t\right)f\left(\vartheta
  ,t\right)^2}{\varepsilon ^2\sigma \left(t\right)^2}{\rm d}X_t,
\end{align*}
where we denoted
\begin{align*}
 q_\varepsilon \left(\vartheta ,t\right)=\frac{\gamma \left(\vartheta
   ,t\right)f\left(\vartheta 
  ,t\right)^2}{\varepsilon ^2\sigma \left(t\right)^2}-a\left(\vartheta
 ,t\right)-\frac{f'\left(\vartheta ,t\right)}{f\left(\vartheta ,t\right)} .
\end{align*}
The initial value is $M\left(\vartheta ,\tau \right)=f\left(\vartheta ,\tau
\right)m\left(\vartheta ,\tau \right)$.

 To describe the asymptotic behavior of $ M(\vartheta ,t)$ and $\dot
 M(\vartheta ,t)$ we need some 
 results from \cite{Kut19}, which we recall here. The function $\gamma
 _*\left(\vartheta ,t\right)=\varepsilon ^{-1}\gamma \left(\vartheta
 ,t\right)$ is bounded and for any $t_0>0$ uniformly on $t\in
 \left[t_0,T\right]$ converges to the function $\gamma _0\left(\vartheta
 ,t\right)=b\left(\vartheta ,t\right)\sigma \left(t\right)f\left(\vartheta
 ,t\right)^{-1}$ (\cite{Kut19}, Lemma 2). It was shown that we can replace the function $q_\varepsilon
 \left(\vartheta ,t\right)$ by the function ${q_*\left(\vartheta
  ,t\right)}/{\varepsilon } $, where
$$
 q_*\left(\vartheta,t\right)= \frac{\gamma _0\left(\vartheta
   ,t\right)f\left(\vartheta ,t\right)^2}{\sigma \left(t\right)^2}=
 \frac{b\left(\vartheta ,t\right)f\left(\vartheta ,t\right)}{\sigma
   \left(t\right)}.
$$
We rewrite the equation  for $M\left(\vartheta ,t\right)$ without changing
notation 
\begin{align*}
{\rm d}M\left(\vartheta ,t\right)&=- \frac{b\left(\vartheta,t
  \right)f\left(\vartheta ,t\right)}{\varepsilon \sigma
  \left(t\right)}M\left(\vartheta ,t\right){\rm d}t+\frac{b\left(\vartheta,t
  \right)f\left(\vartheta ,t\right)}{\varepsilon \sigma \left(t\right)}{\rm
  d}X_t\\ 
&=-\frac{1}{\varepsilon }q_*\left(\vartheta,t\right)M\left(\vartheta
,t\right){\rm d}t+\frac{1}{\varepsilon }q_*\left(\vartheta,t\right){\rm
  d}X_t\\
 &=-\frac{1}{\varepsilon
}q_*\left(\vartheta,t\right)\left[M\left(\vartheta
  ,t\right)-M\left(\vartheta_0 ,t\right)\right]{\rm
  d}t+q_*\left(\vartheta,t\right)\sigma \left(t\right){\rm d}\bar W_t.
\end{align*}
 The equation for derivative is
\begin{align*}
{\rm d}\dot M\left(\vartheta ,t\right)&=-\frac{1}{\varepsilon
}q_*\left(\vartheta,t\right)\dot M\left(\vartheta ,t\right){\rm
  d}t-\frac{1}{\varepsilon }\dot q_*\left(\vartheta,t\right) M\left(\vartheta
,t\right){\rm d}t +\frac{1}{\varepsilon }\dot q_*\left(\vartheta,t\right){\rm
  d}X_t\\
 &=-\frac{1}{\varepsilon }q_*\left(\vartheta,t\right)\dot
M\left(\vartheta ,t\right){\rm d}t+\frac{1}{\varepsilon }\dot
q_*\left(\vartheta,t\right) \left[M\left(\vartheta_0
  ,t\right)-M\left(\vartheta ,t\right)\right]{\rm d}t +\dot
q_*\left(\vartheta,t\right)\sigma \left(t\right) {\rm d}\bar W_t.
\end{align*}

As the contribution of the initial value in $\dot M\left(\vartheta ,t\right) $ is
exponentially small to simplify the exposition we will omit it below. We have
\begin{align}
\label{23-63}
\dot M\left(\vartheta ,t\right)&=\frac{1}{\varepsilon }\int_{\tau }^{t}\phi \left(\vartheta ,s
,t\right)\dot q_*\left(\vartheta,s\right)\left[M\left(\vartheta_0 ,s\right)-M\left(\vartheta
  ,s\right)\right]{\rm d}s \nonumber\\
&\qquad \qquad\qquad \qquad +\int_{\tau }^{t}\phi \left(\vartheta ,s
,t\right)\dot q_* \left(\vartheta ,s\right)\sigma \left(s\right){\rm d}\bar W_s,
\end{align}
where
\begin{align*}
\phi_*\left(\vartheta ,s ,t\right)=\exp\left(-\varepsilon ^{-1}\int_{s}^{t}
q_* \left(\vartheta ,v\right){\rm d}v\right). 
\end{align*}

The asymptotics of ${\dot M\left(\vartheta ,t\right)}$ is described with the
help of the  Lemma 1 in \cite{Kut19}.

Note that
\begin{align*}
\dot M\left(\vartheta _0,t\right)=\int_{\tau }^{t}\phi_*(\vartheta_0,s,t)\dot
q_* (\vartheta_0 ,s)\sigma \left(s\right){\rm d}\bar W_s,
\end{align*}
and
 the asymptotic of stochastic integral  by the same Lemma is
\begin{align*}
&\int_{\tau }^{t}\phi_*(\vartheta,s,t)\dot q_* (\vartheta
,s)\sigma \left(s\right){\rm d}\bar W_s\\
&\qquad  = \sqrt{\varepsilon
}\sqrt{{{b(\vartheta,t)\sigma \left(t\right)}
}{{f(\vartheta,t)}}}\;\left(\frac{\dot
  b(\vartheta,t)}{b(\vartheta,t)}+\frac{\dot
  f(\vartheta,t)}{f(\vartheta,t)} \right)\xi _{t,\varepsilon } 
\left(1+O\left(\sqrt{{\varepsilon }}\right)\right).
\end{align*}
Hence
\begin{align*}
\varepsilon ^{-1/2}\dot M\left(\vartheta
_0,t\right)&=\sqrt{{{b(\vartheta_0,t)\sigma \left(t\right)} 
}{{f(\vartheta_0,t)}}}\;\left(\frac{\dot
  b(\vartheta_0,t)}{b(\vartheta_0,t)}+\frac{\dot
  f(\vartheta_0,t)}{f(\vartheta_0,t)} \right)\xi _{t,\varepsilon } 
\left(1+O\left(\sqrt{{\varepsilon }}\right)\right)\\
&\Longrightarrow \sqrt{{{b(\vartheta_0,t)\sigma \left(t\right)}
}{{f(\vartheta_0,t)}}}\;\left(\frac{\dot
  b(\vartheta_0,t)}{b(\vartheta_0,t)}+\frac{\dot
  f(\vartheta_0,t)}{f(\vartheta_0,t)} \right)\xi _{t } .
\end{align*}

The first integral in \eqref{23-63} we write as follows
\begin{align*}
&\frac{1}{\varepsilon }\int_{\tau }^{t}\phi \left(\vartheta ,s
,t\right)\dot q_*\left(\vartheta,s\right)\left[M\left(\vartheta_0 ,s\right)-M\left(\vartheta
  ,s\right)\right]{\rm d}s\\
&\qquad \qquad  =\frac{\left(\vartheta _0-\vartheta \right)}{\varepsilon }\int_{\tau }^{t}\phi \left(\vartheta ,s
,t\right)\dot q_*\left(\vartheta,s\right) \dot M(\tilde\vartheta ,s){\rm d}s\\
&\qquad \qquad  =\frac{\left(\vartheta _0-\vartheta \right)}{\varepsilon }\int_{\tau }^{t}\phi
  \left(\vartheta ,s 
,t\right)\dot q_*\left(\vartheta,s\right) \dot M(\vartheta_0 ,s){\rm d}s +
  \left(\vartheta _0-\vartheta \right)^2  O\left(1\right)\\
&\qquad \qquad = \sqrt{\varepsilon
  }\left(\vartheta _0-\vartheta \right)O\left(1\right)+ \left(\vartheta _0-\vartheta \right)^2  O\left(1\right).
\end{align*}
where $\left|\tilde\vartheta -\vartheta _0\right|\leq \left|\vartheta
-\vartheta _0\right|$.

We have 
\begin{align}
&\Pb_{\vartheta _0} \left(\sup_{t_0\leq t\leq T}\left|\vartheta
  ^\star_{t,\varepsilon }-\vartheta _0\right|\geq \nu \right)\leq
  \Pb_{\vartheta _0} \left(\left|\check\vartheta _{\tau ,\varepsilon
  }-\vartheta _0\right|\geq \frac{\nu }{3}\right)\nonumber\\ 
&\qquad \qquad
  +\Pb_{\vartheta _0}\left(\sup_{t_0\leq t\leq T}\left|\frac{1}{{\rm I}_\tau^t
    (\check\vartheta _{\tau,\varepsilon } )} \int_{\tau }^{t}\frac{\dot
    M(\check\vartheta _{\tau,\varepsilon },s)}{\sigma \left(s\right)} {\rm
    d}\bar W_s \right|\geq \frac{\nu }{3} \right)\nonumber\\ 
&\qquad\qquad
  +\Pb_{\vartheta _0} \left(\frac{\left|\check\vartheta
    _{\tau ,\varepsilon }-\vartheta _0 \right|}{{\rm I}_\tau^{t_0}
    (\check\vartheta _{\tau,\varepsilon } )} \int_{\tau }^{T}\left|\frac{\dot
    M(\check\vartheta _{\tau,\varepsilon },s)\dot M(\tilde\vartheta
    _{\tau,\varepsilon },s )}{\varepsilon   \sigma \left(s\right)^2}
  \right|{\rm d}s\geq \frac{\nu }{3} \right) .
\label{23-66}
\end{align}
Let us  denote $c_0=\inf_{\vartheta \in\Theta }{\rm I}^{t_0}_{\tau} (\vartheta
) $. Then 
\begin{align*}
&\Pb_{\vartheta _0}\left(\sup_{t_0\leq t\leq T}\left|\frac{1}{{\rm I}_\tau^t
    (\check\vartheta _{\tau,\varepsilon } )} \int_{\tau }^{t}\frac{\dot
    M(\check\vartheta _{\tau,\varepsilon },s)}{ \sigma \left(s\right)} {\rm
    d}\bar W_s \right|\geq \frac{\nu }{3} \right)\\ 
&\qquad \qquad\qquad \qquad \leq
  \Pb_{\vartheta _0}\left(\sup_{t_0\leq t\leq T}\left| \int_{\tau
  }^{t}\frac{\dot M(\check\vartheta _{\tau,\varepsilon },s)}{\sigma
    \left(s\right)} {\rm d}\bar W_s \right|\geq \frac{\nu c_0 }{3}
  \right)\\ 
&\qquad \qquad\qquad \qquad \leq \frac{9}{\nu^2 c_0^2 }\Ex_{\vartheta
    _0} \int_{\tau }^{T}\frac{\dot M(\check\vartheta _{\tau,\varepsilon
    },s)^2}{ \sigma \left(s\right)^2} {\rm d}s\leq C\varepsilon .
\end{align*}
The last integral in \eqref{23-66} is a  bounded in probability function and
the corresponding 
probability tends to zero thanks to the consistency of the estimator
$\check\vartheta _{\tau,\varepsilon } $. Therefore the estimator $\vartheta
  ^\star_{t,\varepsilon } $ is uniformly on $t\in \left[t_0,T\right]$
consistent.

Further
\begin{align*}
\frac{1}{{\rm I}_\tau^t (\check\vartheta _{\tau,\varepsilon } )}& \int_{\tau
}^{t}\frac{\dot M(\check\vartheta _{\tau,\varepsilon },s)}{\varepsilon^{3/2}
  \sigma \left(s\right)^2} \left[M\left(\vartheta _0,s\right)-M(\check\vartheta _{\tau,\varepsilon },s )\right]{\rm
  d}s\\
 &=\frac{\left(\vartheta _0-\check\vartheta _{\tau ,\varepsilon
  }\right)}{{\rm I}_\tau^t (\check\vartheta _{\tau,\varepsilon } )} \int_{\tau
}^{t}\frac{\dot M(\check\vartheta _{\tau,\varepsilon },s)^2}{\varepsilon^{3/2}
  \sigma \left(s\right)^2} {\rm d}s\left( 1+O\left(\varepsilon
^{1/3}\right)\right)\\ 
&=-\frac{\left(\check\vartheta _{\tau ,\varepsilon
  }-\vartheta _0\right)}{\sqrt{\varepsilon}\;{\rm I}_\tau^t
  (\vartheta _0 ) } \int_{\tau }^{t}\frac{ S(\vartheta _0 ,s)
}{ \sigma \left(s\right)}{ }\left[\frac{\dot f(\vartheta _0  ,s)}{f(\vartheta _0 
    ,s)}+\frac{\dot b(\vartheta _0 
    ,s)}{b(\vartheta _0 ,s)} \right]^2\xi
_{s,\varepsilon }^2 {\rm d}s\left( 1+O\left(\varepsilon ^{1/3}\right)\right).
\end{align*}
Recall that (\cite{Kut19}, Lemma 4)
\begin{align}
&\int_{\tau }^{t}\frac{ S(\vartheta _0 ,s)
}{ \sigma \left(s\right)}{ }\left[\frac{\dot f(\vartheta _0  ,s)}{f(\vartheta _0 
    ,s)}+\frac{\dot b(\vartheta _0 
    ,s)}{b(\vartheta _0 ,s)} \right]^2\xi
_{s,\varepsilon }^2 {\rm d}s\nonumber\\
&\qquad \qquad \longrightarrow \int_{\tau }^{t}\frac{ S(\vartheta _0 ,s)
}{2 \sigma \left(s\right)}{ }\left[\frac{\dot f(\vartheta _0  ,s)}{f(\vartheta _0 
    ,s)}+\frac{\dot b(\vartheta _0 
    ,s)}{b(\vartheta _0 ,s)} \right]^2{\rm d}s={\rm I}_\tau^t
  (\vartheta _0 ).
\label{23-67}
\end{align}
Therefore
\begin{align}
\label{23-68}
&\frac{\left(\check\vartheta _{\tau,\varepsilon }-\vartheta _0
    \right)}{\sqrt{\varepsilon }} +\frac{1}{{\rm I}_\tau^t (\check\vartheta
    _{\tau,\varepsilon } )} \int_{\tau }^{t}\frac{\dot M(\check\vartheta
    _{\tau,\varepsilon },s)}{\varepsilon^{3/2} \sigma \left(s\right)^2}
  \left[M\left(\vartheta _0,s\right)-M(\check\vartheta _{\tau,\varepsilon },s )\right]{\rm d}s\nonumber\\
&\quad  =\frac{\left(\check\vartheta _{\tau,\varepsilon }-\vartheta _0
    \right)}{\sqrt{\varepsilon }}\left[1-\frac{1}{{\rm I}_\tau^t (\vartheta
    _0 )}    \int_{\tau }^{t}\frac{ S(\vartheta _0 ,s)
}{ \sigma \left(s\right)}{ }\left[\frac{\dot f(\vartheta _0  ,s)}{f(\vartheta _0 
    ,s)}+\frac{\dot b(\vartheta _0 
    ,s)}{b(\vartheta _0 ,s)} \right]^2\xi
_{s,\varepsilon }^2 {\rm d}s\left( 1+O\left(\varepsilon
    ^{1/3}\right)\right)\right]\nonumber\\
&\quad  =\frac{\left(\check\vartheta _{\tau,\varepsilon }-\vartheta _0
    \right)}{\sqrt{\varepsilon }}O\left(\varepsilon
    ^{1/3}\right)=\frac{O\left(\varepsilon
    ^{2/3}\right)}{\sqrt{\varepsilon }}\longrightarrow 0.
\end{align}
The convergence \eqref{23-67} is sufficient for the asymptotic normality of
the stochastic integral
\begin{align*}
\int_{\tau }^{t}\frac{\dot M(\check\vartheta _{\tau,\varepsilon },s)}{\sqrt{\varepsilon } \sigma
  \left(s\right)} {\rm d}\bar W_s= \int_{\tau }^{t}\frac{\dot M(\vartheta
  _0,s)}{\sqrt{\varepsilon } \sigma \left(s\right)} {\rm d}\bar W_s\left(1+O\left(\varepsilon
^{1/3}\right)\right)\Longrightarrow {\cal N}\left(0,{\rm I}_\tau^t (\vartheta
    _0 )\right).
\end{align*}

Finally we obtained the asymptotic normality of One-step MLE-process: for any
$t>\tau $ 
\begin{align*}
\frac{\vartheta ^\star_{t,\varepsilon }-\vartheta _0}{\sqrt{\varepsilon
}}\Longrightarrow \eta _t\left(\vartheta _0\right)\sim {\cal N}\left(0,{\rm I}_\tau^t (\vartheta
    _0 )^{-1}\right).
\end{align*}

To prove the convergence of moments we slightly modify the calculations in
\eqref{23-68} :
\begin{align*}
&\frac{\left(\check\vartheta _{\tau,\varepsilon }-\vartheta _0
    \right)}{\sqrt{\varepsilon }} +\frac{1}{{\rm I}_\tau^t (\check\vartheta
    _{\tau,\varepsilon } )} \int_{\tau }^{t}\frac{\dot M(\check\vartheta
    _{\tau,\varepsilon },s)}{\varepsilon^{3/2} \sigma \left(s\right)^2}
  \left[M\left(\vartheta _0,s\right)-M(\check\vartheta _{\tau,\varepsilon },s )\right]{\rm d}s\\
&\quad  =\frac{\left(\check\vartheta _{\tau,\varepsilon }-\vartheta _0
    \right)}{\sqrt{\varepsilon }}\left[1-\frac{1}{{\rm I}_\tau^t (\check\vartheta
    _{\tau,\varepsilon } )}    \int_{\tau }^{t}\frac{ S(\check\vartheta
    _{\tau,\varepsilon } ,s)
}{ \sigma \left(s\right)}{ }\left[\frac{\dot f(\check\vartheta
    _{\tau,\varepsilon }  ,s)}{f(\check\vartheta
    _{\tau,\varepsilon } 
    ,s)}+\frac{\dot b(\check\vartheta
    _{\tau,\varepsilon } 
    ,s)}{b(\check\vartheta
    _{\tau,\varepsilon }  ,s)} \right]^2 {\rm d}s\left( 1+O\left(\varepsilon
    ^{1/3}\right)\right)\right]\\
&\qquad + \frac{\left(\check\vartheta _{\tau,\varepsilon }-\vartheta _0
    \right)}{{\rm I}_\tau^t (\check\vartheta
    _{\tau,\varepsilon } )\sqrt{\varepsilon }}   \int_{\tau }^{t}\frac{ S(\check\vartheta
    _{\tau,\varepsilon } ,s)
}{ \sigma \left(s\right)}{ }\left[\frac{\dot f(\check\vartheta
    _{\tau,\varepsilon }  ,s)}{f(\check\vartheta
    _{\tau,\varepsilon } 
    ,s)}+\frac{\dot b(\check\vartheta
    _{\tau,\varepsilon } 
    ,s)}{b(\check\vartheta
    _{\tau,\varepsilon }  ,s)} \right]^2\left\{\xi
_{s,\varepsilon }^2-\frac{1}{2}\right\} {\rm d}s\left( 1+O\left(\varepsilon
    ^{1/3}\right)\right)\\
&\quad =\frac{\left(\check\vartheta _{\tau,\varepsilon }-\vartheta _0
    \right)^2}{\sqrt{\varepsilon }}O\left(1\right) +\left(\check\vartheta _{\tau,\varepsilon }-\vartheta _0
    \right)O\left(1\right).
\end{align*} 
The  terms  $\varepsilon ^{-1/2}\left(\check\vartheta _{\tau,\varepsilon }-\vartheta _0
    \right)^2$,   $O\left(1\right) $ and the integral containing $\xi
_{s,\varepsilon }^2-\frac{1}{2}$ have bounded all polynomial moments
    (see \eqref{23-56} here and  the estimate  (29) in
    \cite{Kut19}). Therefore for any $p>1$ there exists a constant $C>0$ such
    that 
\begin{align}
\label{23-69}
\Ex_{\vartheta _0}\left| \frac{\vartheta^\star
    _{t,\varepsilon } -\vartheta _0}{\sqrt{\varepsilon }}\right|^p\leq C.
\end{align}

The presented here proof can be used for verification of the convergence of
the vectors: for any $K=1,2,\ldots$ and $t_0<t_1<\ldots<t_K\leq T$
\begin{align*}
\left(\eta _{t_1,\varepsilon },\ldots,\eta _{t_K,\varepsilon }\right)\Longrightarrow \left(\eta _{t_1 }\left(\vartheta _0\right),\ldots,\eta _{t_K }\left(\vartheta _0\right)\right).
\end{align*}
It is sufficient to study the convergence of the sum $\sum_{k=1}^{K}\lambda
_k\eta _{t_k,\varepsilon } $.
The estimate
\begin{align*}
\Ex_{\vartheta _0}\left|\eta _{t_1,\varepsilon }-\eta _{t_2,\varepsilon
}\right|^4\leq C_1\left|t_1-t_2\right|^4+C_2\left|t_1-t_2\right|^2 
\end{align*}
with some constants $C_1>0$ and $C_2>0$ can be verified following the same steps
as it was done in \cite{Kut17}, Theorem 1.

Having these properties of $\eta _{t,\varepsilon} $  we
obtain the weak convergence of the random process $\eta _{t,\varepsilon },
t_0\leq t\leq T $ to the random process $\eta _{t }\left(\vartheta _0\right),
t_0\leq t\leq T $.

\end{proof}
\begin{remark}
\label{R23-3}
{\rm 
Note that if $\inf_{\vartheta \in\Theta } |\dot S\left(\vartheta ,\tau 
\right)|>0 $, then the condition ${\scr C}_1$ is fulfilled for any $t_0>\tau $.
}
\end{remark}
\begin{remark}
\label{R23-4}
{\rm The proposed above construction of the One-step MLE $\vartheta
  _{T,\varepsilon }^\star$ and the One-step MLE-process $\vartheta
  _{t,\varepsilon }^\star, \tau <t\leq T$ are not asymptotically efficient
  because the limit variance of the MLE $\hat\vartheta _\varepsilon $ is ${\rm
    I}_0^T\left(\vartheta _0\right)^{-1}$ and the limit variance of $\vartheta
  _{T,\varepsilon }^\star$ is ${\rm I}_\tau ^T\left(\vartheta
  _0\right)^{-1}$. Choosing   sufficiently small $\tau $ the ratio
\begin{align*}
\frac{{\rm I}_0 ^T\left(\vartheta _0\right)}{{\rm I}_\tau ^T\left(\vartheta
_0\right)}\geq 1
\end{align*}
can be made close to 1.

Another possibility is to consider the sequence of problems with $\tau=\tau
_\varepsilon \rightarrow 0$ but slowly. The main condition on the rate of
convergence of a preliminary estimator $\bar\vartheta _{\tau_\varepsilon  ,\varepsilon }$  is 
\begin{align*}
\frac{\left(\bar\vartheta _{\tau_\varepsilon ,\varepsilon }-\vartheta _0\right)^2}{\sqrt{\varepsilon }}\longrightarrow 0.
\end{align*}
The MME $\check\vartheta _{\tau_\varepsilon ,\varepsilon } $ satisfies to this condition 
\begin{align*}
\frac{\left(\check\vartheta _{\tau ,\varepsilon }-\vartheta
  _0\right)^2}{\sqrt{\varepsilon }}\approx  \frac{\varepsilon ^{2/3}}{\tau \varepsilon ^{1/2}}=\frac{\varepsilon ^{1/6}}{\tau }\longrightarrow 0.
\end{align*}
Therefore, if we take $\tau =\tau _\varepsilon =\varepsilon ^{1/12}$, then
\begin{align*}
\frac{\left(\check\vartheta _{\tau_\varepsilon  ,\varepsilon }-\vartheta
  _0\right)^2}{\sqrt{\varepsilon }}\approx {\varepsilon ^{1/12}} \longrightarrow 0
\end{align*}
and we have the asymptotic normality of the corresponding One-step estimator:
for any $t>0$
\begin{align}
\label{23-70}
\frac{\vartheta _{t,\varepsilon }^\star-\vartheta _0}{\sqrt{\varepsilon
}}\Longrightarrow  \hat \eta_t\left(\vartheta _0\right)=\frac{1 }{{\rm I}_0^t
  (\vartheta _{0  } )}\int_{0
}^{t}\frac{\dot S\left(\vartheta _0,s\right)}{\sqrt{2S\left(\vartheta
  _0,s\right)\sigma \left(s\right)}}{\rm d}w\left(s \right)\sim {\cal
  N}\left(0,{\rm I}_0^t
  (\vartheta _{0  } )^{-1} \right).
\end{align}

}
\end{remark}

\section{Adaptive filtration}.  

Consider the partially observed two-dimensional linear stochastic system
\eqref{23-56}, where $\vartheta $ is unknown parameter. The
main problem is the estimation of the hidden component $Y_t$ by observations
$X^t=\left(X_s,0\leq s\leq t\right)$.  We can not use the equations
\eqref{23-57}, \eqref{23-58} directly because the true value of $\vartheta$ is
unknown. It is quite natural to replace in \eqref{23-57}, \eqref{23-58} by
some estimator $\bar\vartheta _\varepsilon $ and the corresponding solution
$m\left(\bar\vartheta _\varepsilon ,t\right) $ of \eqref{23-57} can be an
approximation of the conditional expectation $m\left(\vartheta ,t\right)$. The
properties of the
estimator $\bar\vartheta _\varepsilon $, which we need for such approximation
are: for each $t$ it depends of the observations $X^t$, it can be easily
calculated and is  asymptotically efficient. All these properties has One-step
MLE-process $\vartheta _{t,\varepsilon }^\star, \tau _\varepsilon <t\leq T$
and below we realize this construction and evaluate the error of
approximation.

Let us recall the  notations:  $\tau _\varepsilon =\varepsilon ^{1/12}$,
$t_{i+1}=t_i+\varepsilon ^{1/3}$, $N_{\tau _\varepsilon ,\varepsilon
  }=\left[\varepsilon ^{-1/4}\right]$,
\begin{align}
\Psi _{\tau _\varepsilon }\left(\vartheta \right)&=\int_{0}^{\tau _\varepsilon }f\left(\vartheta
,t\right)^2b\left(\vartheta ,t\right)^2{\rm d}t,\quad  \hat\Psi _{\tau
  _\varepsilon ,\varepsilon }=\sum_{i=0}^{N_{\tau _\varepsilon ,\varepsilon
  }-1}\left(\frac{X_{t_{i+1}+\varepsilon }-X_{t_{i+1}}-X_{t_{i}+\varepsilon
  }+X_{t_{i}}}{\varepsilon } \right)^2,\nonumber\\
\check\vartheta _{\tau _\varepsilon }&=\Psi _{\tau _\varepsilon }^{-1}\left(\hat\Psi _{\tau
  _\varepsilon ,\varepsilon } \right),\qquad \hat\eta _t\left(\vartheta _0\right)={\rm
  I}_0^t\left(\vartheta _0\right)^{-1}\int_{0}^{t}\frac{\dot S\left(\vartheta
  _0,s\right)}{\sqrt{2S\left(\vartheta _0,s\right)\sigma \left(s\right)}}\;{\rm
  d}w\left(s\right),\nonumber\\
 \hat\eta _{t,\varepsilon }\left(\vartheta _0\right)&={\rm
  I}_0^t\left(\vartheta _0\right)^{-1}\int_{0}^{t}\frac{\dot S\left(\vartheta
  _0,s\right)}{\sqrt{S\left(\vartheta _0,s\right)\sigma
     \left(s\right)}}\;\xi^{\left(1\right)} _{s,\varepsilon }\;{\rm  d}\bar
 W_s,\qquad \xi^{\left(1\right)} _{t,\varepsilon
 }=\sqrt{\frac{q_t}{\varepsilon }}\int_{0}^{\frac{tq_t}{\varepsilon }}e^{-v} {\rm d}\bar W_{t-\frac{v\varepsilon }{q_t}},\nonumber\\ 
\label{23-71}
{\rm d}M\left(\vartheta ,t\right)&=- \frac{b\left(\vartheta,t
  \right)f\left(\vartheta ,t\right)}{\varepsilon \sigma
  \left(t\right)}M\left(\vartheta ,t\right){\rm d}t+\frac{b\left(\vartheta,t
  \right)f\left(\vartheta ,t\right)}{\varepsilon \sigma \left(t\right)}{\rm
  d}X_t,\\
\label{23-72}
{\rm d}\dot M\left(\vartheta ,t\right)&=- \frac{b\left(\vartheta
 ,t \right)f\left(\vartheta ,t\right)}{\varepsilon \sigma
  \left(t\right)}\dot M\left(\vartheta ,t\right){\rm d}t\nonumber\\
&\qquad \qquad +\frac{\dot b\left(\vartheta,t
  \right)f\left(\vartheta ,t\right)+\dot f\left(\vartheta ,t\right)b\left(\vartheta,t
  \right)}{\varepsilon \sigma \left(t\right)}\left[{\rm
  d}X_t-M\left(\vartheta ,t\right){\rm d}t\right],\\
\vartheta _{t,\varepsilon }^\star&=\check\vartheta _{\tau _\varepsilon }+{\rm
  I}_{\tau _\varepsilon }^t\left(\check\vartheta _{\tau _\varepsilon
}\right)^{-1}\int_{\tau _\varepsilon }^{t}\frac{\dot M\left(\check\vartheta
  _{\tau _\varepsilon },s\right)}{\varepsilon \sigma \left(s\right)^2}\left[{\rm 
  d}X_s-M\left(\check\vartheta _{\tau _\varepsilon },s\right){\rm
    d}s\right],\quad \tau _\varepsilon <t\leq T.
\label{23-73}
\end{align}
For the calculation $\vartheta _{t,\varepsilon }^\star $ we use twice   the
equation \eqref{23-71}. The first time it is considered on the interval
$\left[0,\tau _\varepsilon \right]$ and is rewritten like  \eqref{23-59} to
provide the initial value $M\left(\check\vartheta _{\tau _\varepsilon
  ,\varepsilon },\tau _\varepsilon \right)$ for the solution $M\left(\check\vartheta _{\tau _\varepsilon },t\right)$   of \eqref{23-71}
on the interval $\left[\tau _\varepsilon ,T\right]$, i.e.,
\begin{align*}
{\rm d}M\left(\check\vartheta _{\tau _\varepsilon } ,t\right)&=-
\frac{b\left(\check\vartheta _{\tau _\varepsilon },t
  \right)f\left(\check\vartheta _{\tau _\varepsilon } ,t\right)}{\varepsilon
  \sigma \left(t\right)}M\left(\check\vartheta _{\tau _\varepsilon }
,t\right){\rm d}t+\frac{b\left(\check\vartheta _{\tau _\varepsilon },t
  \right)f\left(\check\vartheta _{\tau _\varepsilon } ,t\right)}{\varepsilon
  \sigma \left(t\right)}{\rm d}X_t,\quad \tau _\varepsilon \leq t\leq T.
\end{align*}
The similar procedure we have for the solution $ \dot M\left(\check\vartheta
_{\tau _\varepsilon } ,t\right)$ of the equation \eqref{23-72}. Having these
two processes the estimator $\vartheta _{t,\varepsilon }^\star $ can be easily 
calculated. 

The {\it adaptive filtration}  for this model of observations is introduced by
the equation for $\hat m_t$ as follows
\begin{align}
\label{23-74}
{\rm d}\hat m_t=-
\frac{b\left(\vartheta _{t,\varepsilon}^\star,t
  \right)f\left(\vartheta _{t,\varepsilon}^\star ,t\right)}{\varepsilon
  \sigma \left(t\right)}\hat m_t{\rm d}t+\frac{b\left(\vartheta _{t,\varepsilon}^\star,t
  \right)}{\varepsilon
  \sigma \left(t\right)}{\rm d}X_t,\quad \tau _\varepsilon \leq t\leq T
\end{align}
with initial value $\hat m_{\tau _\varepsilon }=m\left(\check \vartheta _{\tau
  _\varepsilon },\tau _\varepsilon \right) $. Note that the calculation of
$\hat m_t$ does not require the solution of Riccati equation. 

Let us discuss the errors of approximations: $Y_t-m\left(\vartheta
,t\right)$ and $Y_t-\hat m_t$. Recall that  $Y_t$ is solution of the
equation
\begin{align*}
{\rm d}Y_t=a\left(\vartheta_0 ,t\right)Y_t{\rm d}t+b\left(\vartheta_0 ,t\right){\rm d}V_t,\qquad Y_0=y_0,\quad 0\leq t\leq T.
\end{align*}

Introduce the set
\begin{align*}
\AA=\left(\omega \;:\; \left|\vartheta _{t,\varepsilon}^\star-\vartheta
_0\right|\leq  \varepsilon ^{\frac{3}{8}}\right).
\end{align*}
and the random variables $\xi
  _{t,\varepsilon }^{\left(1\right)},\xi
  _{t,\varepsilon }^{\left(2\right)} $, which are  asymptotically normal
\begin{align*}
\xi _{t,\varepsilon }^{\left(1\right)}\Longrightarrow \xi
  _{t }^{\left(1\right)}\sim {\cal N}\left(0,\frac{1}{2}\right),\qquad \quad \xi _{t,\varepsilon }^{\left(2\right)}\Longrightarrow \xi
  _{t }^{\left(2\right)}\sim {\cal N}\left(0,\frac{1}{2}\right),
\end{align*}
where $\xi
  _{t }^{\left(1\right)},\xi
  _{t }^{\left(2\right)}, 0<t\leq T$  are independent random variables. These
random variables are independent of the Gaussian process $\hat\eta _t\left(\vartheta _0\right),0<t\leq
T$  too. 

Introduce the condition

$\tilde{\scr C}_1$. {\it The functions $f\left(\vartheta ,t\right)$ and
$b\left(\vartheta ,t\right)$ are such that}
\begin{align*}
\inf_{\vartheta \in\Theta }\left| \dot f\left(\vartheta ,0\right)b\left(\vartheta
,0\right)+f\left(\vartheta ,0\right)\dot b\left(\vartheta
,0\right)\right|>0 .
\end{align*}
Note that under this condition we have  $\dot\Psi _\tau \left(\vartheta
\right)>0 $  and  ${\rm I}_0^\tau \left(\vartheta
\right)>0$ for any $\tau >0$.

\begin{theorem}
\label{T23-2} 
Let the conditions ${\scr A},\tilde{\scr C}_1,{\scr C}_2$     be
fulfilled, then, 
\begin{enumerate}
\item  If $\vartheta _0$ is known
then 
\begin{align*}
\frac{Y_t-m\left(\vartheta_0 ,t\right)}{\sqrt{\varepsilon
}}=\sqrt{\frac{b\left(\vartheta_0,t \right)\sigma
    \left(t\right)}{f\left(\vartheta _0,t\right)}}\;\left(\xi _{t,\varepsilon
}^{\left(1\right)}-\xi _{t,\varepsilon }^{\left(2\right)}
\right)\left(1+O\left(\varepsilon ^{1/2}\right)\right).
\end{align*}
\item  If $\vartheta _0$ is unknown, then on the  set $\AA$ we have the
  representation 
\begin{align*}
\frac{Y_t-\hat m_t}{\sqrt{\varepsilon }}=\frac{\dot f\left(\vartheta
  _0,t\right)}{f\left(\vartheta 
  _0,t\right)} \;\hat\eta _{t,\varepsilon }\left(\vartheta _0\right)\,Y_t+\sqrt{\frac{b\left(\vartheta_0,t
    \right)\sigma \left(t\right)}{f\left(\vartheta _0,t\right) }}\;\left(\xi
  _{t,\varepsilon }^{\left(1\right)}-\xi 
  _{t,\varepsilon }^{\left(2\right)} \right)+O\left(\varepsilon ^{q_*}\right)
\end{align*}
where $q_*>0$ and  
\begin{align*}
\Ex_{\vartheta _0}\Bigl(\hat\eta _{t,\varepsilon }\left(\vartheta _0\right)\,Y_t\Bigr)^2\longrightarrow
\Ex_{\vartheta _0}\hat\eta _{t }\left(\vartheta _0\right)^2\, \Ex_{\vartheta
  _0}   Y_t ^2 ={\rm
  I}_{0}^t\left(\vartheta _{0
}\right)^{-1}  \, \Ex_{\vartheta
  _0}   Y_t ^2 .
\end{align*}
For any
$p_*>1$  there exists a constant $C_p>0$ such that 
\begin{align*}
\Pb_{\vartheta _0}\left(\AA^c\right)\leq C_p\,\varepsilon ^{p_*}
\end{align*}
\end{enumerate}
\end{theorem}
\begin{proof}
 Suppose that $\vartheta_0 $ is known and we are interested by the error $Y_t-m\left(\vartheta _0 ,t\right)  $. Of course, we have immediately: for any $t\in (0,T]$
\begin{align*}
\Ex_{\vartheta_0 }\left(Y_t-m\left(\vartheta _0 ,t\right)\right)^2=\gamma
\left(\vartheta _0 ,t\right) =\frac{b\left(\vartheta _0 ,t\right)\sigma
  \left(t\right)}{f\left(\vartheta _0 ,t\right)} \,\varepsilon \left(1+o\left(\varepsilon \right)\right).
\end{align*}
Further, using equation  \eqref{23-57} and asymptotics of
$\gamma \left(\vartheta ,t\right)$, we can write the equation for $\delta \left(t\right)=Y_t-m\left(\vartheta _0,t\right)$
\begin{align*}
{\rm d}\delta \left(t\right)&=a\left(\vartheta_0 ,t\right)\delta \left(t\right){\rm d}t-\frac{\gamma
  \left(\vartheta_0 ,t\right)f\left(\vartheta_0 ,t\right)^2}{\varepsilon ^2\sigma
  \left(t\right)^2}\delta \left(t\right) {\rm d}t+b\left(\vartheta_0 ,t\right){\rm d}V_t-\frac{\gamma
  \left(\vartheta_0 ,t\right)f\left(\vartheta_0 ,t\right)}{\varepsilon \sigma
  \left(t\right)} {\rm d}W_t\\
&=-q_\varepsilon \left(\vartheta_0 ,t\right)\delta \left(t\right){\rm d}t+b\left(\vartheta_0 ,t\right){\rm d}V_t-\frac{\gamma
  \left(\vartheta_0 ,t\right)f\left(\vartheta_0 ,t\right)}{\varepsilon \sigma
  \left(t\right)} {\rm d}W_t\\
&=-\frac{1}{\varepsilon }q_* \left(\vartheta_0 ,t\right)\delta \left(t\right){\rm d}t+b\left(\vartheta_0 ,t\right){\rm d}V_t-b\left(\vartheta_0 ,t\right){\rm d}W_t,
\end{align*}
where we omitted the term $O\left(\varepsilon \right)$.

Hence for any $t>0$ we have the relation
\begin{align}
\label{23-75}
Y_t-m\left(\vartheta_0 ,t\right)=\sqrt{\frac{b\left(\vartheta_0
  ,t\right)\sigma \left(t\right)}{f\left(\vartheta_0 ,t\right)}}\left(\xi
_{t,\varepsilon }^{\left(1\right)}-\xi
_{t,\varepsilon }^{\left(2\right)}   \right) \sqrt{\varepsilon
}\left(1+O\left(\sqrt{\varepsilon }\right)\right).
\end{align}

The equation for $\hat m_t$ is
\begin{align*}
{\rm d}\hat m_t&=-\frac{b\left(\vartheta _{t,\varepsilon}^\star,t \right)f\left(\vartheta
  _{t,\varepsilon}^\star ,t\right)}{\varepsilon \sigma \left(t\right)}\hat
m_t{\rm d}t+\frac{b\left(\vartheta _{t,\varepsilon}^\star,t
  \right)f\left(\vartheta _0,t \right)}{\varepsilon \sigma
  \left(t\right)}Y_t{\rm d}t+b\left(\vartheta
_{t,\varepsilon}^\star,t \right){\rm d}W_t
\end{align*}
Hence
\begin{align*}
{\rm d}\left(Y_t-\hat m_t\right)&= a\left(\vartheta _0,t\right)Y_t{\rm d}t-\frac{b\left(\vartheta _{t,\varepsilon}^\star,t
  \right)\left[f\left(\vartheta _0,t \right)-f\left(\vartheta _{t,\varepsilon}^\star,t
  \right)\right]}{\varepsilon \sigma
  \left(t\right)}Y_t{\rm d}t-b\left(\vartheta
_{t,\varepsilon}^\star,t \right){\rm d}W_t\\
&\qquad  -\frac{b\left(\vartheta _{t,\varepsilon}^\star,t \right)f\left(\vartheta
  _{t,\varepsilon}^\star ,t\right)}{\varepsilon \sigma \left(t\right)}\left(Y_t-\hat
m_t\right){\rm d}t+b\left(\vartheta _0,t\right){\rm d}V_t.
\end{align*}
Further
\begin{align*}
-\varepsilon \ln \phi_*\left(\vartheta ^\star_\varepsilon ,s,t\right)&=\int_{s
}^{t}\frac{S\left(\vartheta _{v,\varepsilon}^\star,v \right)}{\sigma
  \left(v\right)}{\rm d}v=\left(t-s\right)\frac{S\left(\vartheta _0,t
  \right)}{\sigma \left(t\right)}+ \int_{s }^{t}\left(\frac{S\left(\vartheta
  _0,v \right)}{\sigma \left(v\right)}-\frac{S\left(\vartheta _0,t \right)}{\sigma \left(t\right)}
\right){\rm d}v\\
 &\qquad +\int_{s }^{t}\frac{S\left(\vartheta
  _{v,\varepsilon}^\star,v \right)- S\left(\vartheta _0,v \right)}{\sigma
  \left(v\right)}{\rm d}v\\ 
&=\left(t-s\right)q\left(\vartheta_0,t
\right)+\int_{s }^{t}\left(v-t\right)\frac{S'(\vartheta _0,\tilde v
  )\sigma \left(\tilde v\right)-S(\vartheta _0,\tilde v
  )\sigma' \left(\tilde  v\right)}{\sigma \left(\tilde v\right)^2} {\rm d}v\\
&\qquad +\int_{s }^{t} \left(\vartheta
_{v,\varepsilon}^\star-\vartheta _0 \right)\frac{\dot S(\tilde\vartheta_v,v
  )}{\sigma \left(v\right)}{\rm d}v.
\end{align*}
We have elementary estimates
\begin{align*}
\left|\int_{s }^{t}\left(v-t\right)\frac{S'(\vartheta _0,\tilde v
  )\sigma \left(\tilde v\right)-S(\vartheta _0,\tilde v
  )\sigma' \left(\tilde  v\right)}{\sigma \left(v\right)^2} {\rm d}v \right|&\leq C\left(t-s\right)^2,\\
\left|\int_{s }^{t} \left(\vartheta
_{v,\varepsilon}^\star-\vartheta _0 \right)\frac{\dot S(\tilde\vartheta_v,v
  )}{\sigma \left(v\right)}{\rm d}v\right|&\leq C\left(t-s\right)\sup_{s\leq v\leq t}\left|\vartheta
_{v,\varepsilon}^\star-\vartheta _0 \right|.
\end{align*}
Introduce the set 
\begin{align*}
\AA_\varepsilon =\left(\omega : \;\sup_{s\leq v\leq t}\left|\vartheta
_{v,\varepsilon}^\star-\vartheta _0 \right|\leq \varepsilon ^{1/2 -\nu _*}\right),
\end{align*}
where $\nu _*>0$ will be chosen later. On the set $\AA_\varepsilon$  we can write
\begin{align}
\label{23-76}
\left|\phi_*\left(\vartheta
^\star_\varepsilon ,s,t\right) \exp\left(\frac{t-s}{\varepsilon
}q\left(\vartheta _0,t\right)\right)-1\right|\leq C \frac{\left(t-s\right)^2}{\varepsilon
}+C \frac{\left(t-s\right)}{\varepsilon }\;\varepsilon ^{\frac{1}{2}-\nu _*}.
\end{align}

We have the representation
\begin{align*}
Y_t-\hat m_t&=  \frac{1}{\varepsilon }\int_{\tau _\varepsilon  }^{t}\phi_*\left(\vartheta
^\star_\varepsilon ,s,t\right) \left[a\left(\vartheta _0,s\right)\varepsilon +
\frac{b\left(\vartheta _{s,\varepsilon}^\star,s
  \right)\left[f\left(\vartheta _{s,\varepsilon}^\star,s
  \right)-f\left(\vartheta _0,s \right)\right]}{ \sigma
  \left(s\right)}\right]Y_s{\rm d}s\\
&\qquad -\int_{\tau _\varepsilon  }^{t}\phi_*\left(\vartheta
^\star_\varepsilon ,s,t\right)b\left(\vartheta _{s,\varepsilon}^\star,s
  \right){\rm d}W_s+\int_{\tau _\varepsilon  }^{t}\phi_*\left(\vartheta
^\star_\varepsilon ,s,t\right)b\left(\vartheta _0,s
  \right){\rm d}V_s.
\end{align*}

Denote
\begin{align*}
B\left(\vartheta ,\vartheta _0,s\right)&=a\left(\vartheta _0,s\right)\varepsilon +
\frac{b\left(\vartheta ,s
  \right)\left[f\left(\vartheta,s
  \right)-f\left(\vartheta _0,s \right)\right]}{ \sigma
  \left(s\right)},\\
B_*\left(\vartheta 
,\vartheta _0,s\right)&=
\frac{b\left(\vartheta ,s
  \right)\left[f\left(\vartheta,s
  \right)-f\left(\vartheta _0,s \right)\right]}{ \sigma
  \left(s\right)},\qquad \quad \bar \kappa =\frac{\inf_{\vartheta
    ,s}S\left(\vartheta ,s\right)}{\sup_{s}\sigma \left(s\right)}, 
\end{align*}
 and fix some small $\delta >0$. Then we can write the estimate
\begin{align*}
&\frac{1}{\varepsilon }\int_{\tau _\varepsilon }^{t-\delta
}\phi_*\left(\vartheta ^\star_\varepsilon ,s,t\right)\left|B\left(\vartheta
_{s,\varepsilon}^\star,\vartheta _0,s\right)Y_s\right|{\rm d}s\leq
\frac{1}{\varepsilon } \int_{\tau _\varepsilon }^{t-\delta
} \left|B\left(\vartheta
_{s,\varepsilon}^\star,\vartheta _0,s\right)Y_s\right|{\rm
  d}s\;\exp\left(-\frac{\bar\kappa \;\delta }{\varepsilon }\right).
\end{align*}
This allows us to write
\begin{align*}
&\frac{1}{\varepsilon }\int_{\tau _\varepsilon  }^{t}\phi_*\left(\vartheta
^\star_\varepsilon ,s,t\right)B\left(\vartheta
_{s,\varepsilon}^\star,\vartheta _0,s\right)Y_s{\rm d}s=\frac{1}{\varepsilon
  }\int_{t-\delta }^{t }\phi_*\left(\vartheta 
^\star_\varepsilon ,s,t\right)B\left(\vartheta
_{s,\varepsilon}^\star,\vartheta _0,s\right){\rm d}s\,Y_t\\
&\qquad \qquad +\frac{1}{\varepsilon }\int_{t-\delta   }^{t}\phi_*\left(\vartheta
^\star_\varepsilon ,s,t\right)B\left(\vartheta
_{s,\varepsilon}^\star,\vartheta _0,s\right)\left[Y_s-Y_t\right]{\rm d}s
  +O\left( e^{-\frac{\bar\kappa \delta }{2\varepsilon }}\right).
\end{align*}

On the set $\AA_\varepsilon$ for the first integral in the RHS according to
\eqref{23-76}  we have  
\begin{align*}
&I_{1,\varepsilon }=\frac{1}{\varepsilon }\int_{t-\delta }^{t }\phi_*\left(\vartheta
^\star_\varepsilon ,s,t\right)B\left(\vartheta
_{s,\varepsilon}^\star,\vartheta _0,s\right){\rm d}s\\
&\qquad =\frac{1}{\varepsilon }\int_{t-\delta }^{t
  }\exp\left(-\frac{t-s}{\varepsilon }q_*\left(\vartheta_0,t
  \right)\right)\left(1+O\left(\frac{\delta ^2}{\varepsilon
  }\right)+O\left(\frac{\delta  }{\varepsilon^{\frac{1}{2}+\nu
      _*}}\right)\right)B\left(\vartheta  
_{s,\varepsilon}^\star,\vartheta _0,s\right){\rm d}s\\
&\qquad =\frac{1}{\varepsilon }\int_{t-\delta }^{t
  }\exp\left(-\frac{t-s}{\varepsilon }q_*\left(\vartheta
  _0,t\right)\right)B\left(\vartheta  
_{s,\varepsilon}^\star,\vartheta _0,s\right){\rm d}s+O\left(\frac{\delta ^2}{\varepsilon 
  }\right)+O\left(\frac{\delta}{ \varepsilon^{\frac{1}{2}+\nu
      _*} }\right).
\end{align*}
Further
\begin{align*}
&\frac{1}{\varepsilon }\int_{t-\delta }^{t }\exp\left(-\frac{t-s}{\varepsilon
}q_*\left(\vartheta _0,t\right)\right)B\left(\vartheta
_{s,\varepsilon}^\star,\vartheta _0,s\right){\rm d}s\\
& \qquad  \quad =\frac{1}{\varepsilon }\int_{t-\delta }^{t
  }\exp\left(-\frac{t-s}{\varepsilon 
}q_*\left(\vartheta _0,t\right)\right)B\left(\vartheta
_0,\vartheta _0,s\right){\rm d}s\\
&\quad\qquad  \qquad +\frac{1}{\varepsilon }\int_{t-\delta }^{t
  }\exp\left(-\frac{t-s}{\varepsilon 
}q_*\left(\vartheta _0,t\right)\right)\left(\vartheta
_{s,\varepsilon}^\star-\vartheta _0\right) \dot  B(\tilde\vartheta
  _s,\vartheta _0,s){\rm d}s\\ 
& \qquad  \quad =\frac{\dot B_*\left(\vartheta
_0,\vartheta _0,t\right)\sigma \left(t\right)}{S\left(\vartheta
    _0,t\right)}\;\eta_{t,\varepsilon }\left(\vartheta
  _0\right)\,\sqrt{\varepsilon }+O\left(\varepsilon 
  \right) =\frac{\dot f\left(\vartheta 
_0,t\right)   }{f\left(\vartheta
    _0,t\right)}\;\eta_{t,\varepsilon }\left(\vartheta
  _0\right)\,\sqrt{\varepsilon }+O\left(\varepsilon \right) . 
\end{align*}
Recall that $B(\vartheta _0,\vartheta _0,s)=0$.
The second integral is estimated as follows
\begin{align*}
&\frac{1}{\varepsilon }\int_{t-\delta   }^{t}\phi_*\left(\vartheta
^\star_\varepsilon ,s,t\right)B\left(\vartheta
_{s,\varepsilon}^\star,\vartheta _0,s\right)\left[Y_s-Y_t\right]{\rm d}s\\
&\qquad \qquad =\frac{1}{\varepsilon }\int_{t-\delta }^{t }\exp\left(-\frac{t-s}{\varepsilon
}q_*\left(\vartheta _0,t\right)\right)B\left(\vartheta
_{s,\varepsilon}^\star,\vartheta _0,s\right)\int_{s}^{t}b\left(\vartheta
  _0,r\right){\rm d}V_r\;{\rm d}s +O\left(\delta \right)\\
&\qquad \qquad =\frac{1}{\varepsilon }\int_{t-\delta }^{t }\exp\left(-\frac{t-s}{\varepsilon
}q_*\left(\vartheta _0,t\right)\right)B\left(\vartheta
_{s,\varepsilon}^\star,\vartheta _0,s\right)b\left(\vartheta
  _0,t\right)\frac{\left(V_t-V_s\right)}{\sqrt{\delta }}\;{\rm d}s\sqrt{\delta
  }\; +O\left(\delta \right)\\ 
&\qquad \qquad =O\left(\sqrt{\varepsilon \delta }\right) +O\left(\delta \right).
\end{align*}
Here we used the relations $B\left(\vartheta
_{s}^\star,\vartheta _0,s\right)=\left(\bar \vartheta
_{s,\varepsilon}^\star-\vartheta _0\right)  \dot B\left(\bar\vartheta
,\vartheta _0,s\right)=\dot B\left(\bar\vartheta ,\vartheta _0,s\right)\eta
_{s,\varepsilon }\left(\vartheta _0\right)\sqrt{\varepsilon } $. 
For the stochastic integral the similar relations are
\begin{align*}
\int_{\tau _\varepsilon }^{t-\delta
  }\phi_*\left(\vartheta ^\star_\varepsilon ,s,t\right)b\left(\vartheta
  _{s,\varepsilon}^\star,s \right){\rm d}W_s=O\left(e^{-\frac{\bar\kappa \delta }{\varepsilon }}\right)
\end{align*}
and
\begin{align*}
&\int_{t-\delta }^{t}\phi_*\left(\vartheta ^\star_\varepsilon
  ,s,t\right)b\left(\vartheta _{s,\varepsilon}^\star,s \right){\rm
    d}W_s\\
&\qquad   = \int_{t-\delta
  }^{t }\exp\left(-\frac{t-s}{\varepsilon }q_*\left(\vartheta
  _0,t\right)\right)\left(1+O\left(\frac{\delta 
    ^2}{\varepsilon }\right)+O\left(\frac{\delta }{\varepsilon
    ^{\frac{1}{2}+\nu _*}}\right)\right)b\left(\vartheta
  _ {s,\varepsilon}^\star,s \right){\rm
    d}W_s\\ & \qquad = \int_{t-\delta 
  }^{t }\exp\left(-\frac{t-s}{\varepsilon }q_*\left(\vartheta
  _0,t\right)\right) b\left(\vartheta 
  _0,s \right) {\rm
    d}W_s+O\left(\frac{\delta 
    ^2}{\sqrt{\varepsilon }}\right)+O\left(\frac{\delta }{\varepsilon
    ^{\nu _*}}\right)+O\left({\varepsilon^{1-\nu _*} }\right)\\
 & \qquad =\sqrt{\frac{ b\left(\vartheta
  _0,t \right)\sigma \left(t\right)}{{f\left(\vartheta _0,t\right)  }}}\;\xi
  _{t,\varepsilon 
  }^{\left(2\right)}\;\sqrt{\varepsilon } +O\left(\frac{\delta 
    ^2}{\sqrt{\varepsilon }}\right)+O\left(\frac{\delta }{\varepsilon
    ^{\nu _*}}\right)+O\left({\varepsilon^{1-\nu _*} }\right).
\end{align*}
For the second stochastic integral we have the similar presentation
\begin{align*}
&\int_{\tau _\varepsilon }^{t}\phi_*\left(\vartheta ^\star_\varepsilon
  ,s,t\right)b\left(\vartheta _{s,\varepsilon}^\star,s \right){\rm
    d}V_s\\
&\qquad \qquad =\sqrt{\frac{ b\left(\vartheta
  _0,t \right)\sigma \left(t\right)}{{f\left(\vartheta _0,t\right)  }}}\;\xi
  _{t,\varepsilon 
  }^{\left(1\right)}\;\sqrt{\varepsilon } +O\left(\frac{\delta 
    ^2}{\sqrt{\varepsilon }}\right)+O\left(\frac{\delta }{\varepsilon
    ^{\nu _*}}\right)+O\left({\varepsilon^{1-\nu _*} }\right).
\end{align*}

Finally, if we  set $\delta =\varepsilon ^{\frac{7}{8}}$ and $\nu
_*={\frac{1}{8}}$ then  on the set $\AA_\varepsilon=\AA$ we obtain the representation 
\begin{align*}
\frac{Y_t-\hat m_t}{\sqrt{\varepsilon }}&=\frac{\dot f\left(\vartheta
  _0,t\right) }{f\left(\vartheta
  _0,t\right)}\;\eta_{t,\varepsilon }(\vartheta _0)\,Y_t+ \sqrt{\frac{ b\left(\vartheta _0,t
    \right)\sigma \left(t\right)}{{f\left(\vartheta _0,t\right) }}}\;\left(\xi
_{t,\varepsilon }^{\left(2\right)}-\xi _{t,\varepsilon
}^{\left(1\right)}\right) +O\left(\varepsilon ^{\frac{1}{4}}\right).
\end{align*}

Let us estimate the probability $\Pb_{\vartheta _0}\left(\AA^c\right)$. Denote
\begin{align*}
c_1=\inf_{\vartheta\in\Theta  }{\rm
  I}_{\tau _\varepsilon }^{t-\delta }\left(\vartheta \right)
\end{align*}
 We have
\begin{align*}
&\Pb_{\vartheta _0}\left(\sup_{t-\delta \leq s\leq t}\frac{\left|\vartheta
_{s,\varepsilon }^\star-\vartheta 
_0\right|}{\sqrt{\varepsilon }}\geq \varepsilon^{-\frac{1}{12}}  \right) \leq \Pb_{\vartheta _0}\left(\frac{1}{{\rm
  I}_{\tau _\varepsilon }^{t-\delta }\left(\check\vartheta _{\tau _\varepsilon
}\right) } \sup_{t-\delta \leq s\leq t}\left| \int_{\tau _\varepsilon }^{s} \frac{\dot M \left(\check\vartheta _{\tau _\varepsilon
},r\right) }{\sqrt{\varepsilon} \sigma \left(r\right)}{\rm d}\bar W_r \right|\geq \frac{
    \varepsilon^{-\frac{1}{12}}}{2}      \right)\\
&\qquad +
\Pb_{\vartheta _0}\left(\frac{\left|\check\vartheta _{\tau
      _\varepsilon ,\varepsilon }-\vartheta _0\right|}{\sqrt{\varepsilon }\;{\rm
  I}_{\tau _\varepsilon }^{t-\delta }\left(\check\vartheta _{\tau _\varepsilon
}\right) }\sup_{t-\delta \leq s\leq t}\Bigl|{\rm
  I}_{\tau _\varepsilon }^s\left(\check\vartheta _{\tau _\varepsilon
}\right)-\int_{\tau _\varepsilon }^{s} \frac{\dot M \left(\check\vartheta _{\tau _\varepsilon
},r\right) \dot M (\tilde\vartheta _r,r)}{\varepsilon \sigma \left(r\right)^2}{\rm d}r \Bigr|\geq \frac{
    \varepsilon^{-\frac{1}{12}}}{2}  \right).
\end{align*}
To estimate  the first probability we use BGD inequality \cite{BDG72}: for any $p>0$
\begin{align}
\label{23-77}
&\Pb_{\vartheta _0}\left(\frac{1}{{\rm
  I}_{\tau _\varepsilon }^{t-\delta }\left(\check\vartheta _{\tau _\varepsilon
}\right) } \sup_{t-\delta \leq s\leq t}\left| \int_{\tau _\varepsilon }^{s}
  \frac{\dot M \left(\check\vartheta _{\tau _\varepsilon 
},r\right) }{\sqrt{\varepsilon} \sigma \left(r\right)}{\rm d}\bar W_r
  \right|\geq \frac{   \varepsilon^{-\frac{1}{12}}}{2}      \right)\nonumber\\
&\qquad \qquad \leq \Pb_{\vartheta _0}\left(\sup_{t-\delta \leq s\leq t}\left| \int_{\tau _\varepsilon }^{s}
  \frac{\dot M \left(\check\vartheta _{\tau _\varepsilon 
},r\right) }{\sqrt{\varepsilon} \sigma \left(r\right)}{\rm d}\bar W_r
  \right|\geq \frac{ c_1  \varepsilon^{-\frac{1}{12}}}{2}  \right)\nonumber\\
&\qquad \qquad \leq  C\Ex_{\vartheta _0}\left| \int_{\tau _\varepsilon }^{t}
  \frac{\dot M \left(\check\vartheta _{\tau _\varepsilon 
},r\right)^2 }{{\varepsilon} \sigma \left(r\right)^2}{\rm d} r \right|^p\varepsilon ^\frac{p}{6}\leq C\;\varepsilon ^\frac{p}{6}
\end{align}
Using  the relations \eqref{23-68} and  the estimate  (21) in \cite{Kut19} we
obtain the estimate for the second probability 
\begin{align*}
&\Pb_{\vartheta _0}\left(\frac{\left|\check\vartheta _{\tau _\varepsilon
      ,\varepsilon }-\vartheta _0\right|}{\sqrt{\varepsilon }\;{\rm I}_{\tau
      _\varepsilon }^{t-\delta }\left(\check\vartheta _{\tau _\varepsilon
    }\right) }\sup_{t-\delta \leq s\leq t}\Bigl|{\rm I}_{\tau _\varepsilon
  }^s\left(\check\vartheta _{\tau _\varepsilon }\right)-\int_{\tau
    _\varepsilon }^{s} \frac{\dot M \left(\check\vartheta _{\tau _\varepsilon
    },r\right) \dot M (\tilde\vartheta _r,r)}{\varepsilon \sigma
    \left(r\right)^2}{\rm d}r \Bigr|\geq \frac{
    \varepsilon^{-\frac{1}{12}}}{2} \right)\\
 &\qquad \leq \Pb_{\vartheta
    _0}\left(\frac{\left|\check\vartheta _{\tau _\varepsilon ,\varepsilon
    }-\vartheta _0\right|}{\sqrt{\varepsilon } }\sup_{t-\delta \leq s\leq
    t}\Bigl|{\rm I}_{\tau _\varepsilon }^s\left(\check\vartheta _{\tau
    _\varepsilon }\right)-\int_{\tau _\varepsilon }^{s} \frac{\dot M
    \left(\check\vartheta _{\tau _\varepsilon },r\right) \dot M
    (\tilde\vartheta _r,r)}{\varepsilon \sigma \left(r\right)^2}{\rm d}r
  \Bigr|\geq \frac{c_1 \varepsilon^{-\frac{1}{12}}}{2} \right)\\
&\qquad \leq \Pb_{\vartheta
    _0}\left(\frac{\left|\check\vartheta _{\tau _\varepsilon ,\varepsilon
    }-\vartheta _0\right|}{\sqrt{\varepsilon } }\sup_{t-\delta \leq s\leq
    t}\Bigl|{\rm I}_{\tau _\varepsilon }^s\left(\vartheta _0\right)-\int_{\tau _\varepsilon }^{s} \frac{\dot M
    \left(\vartheta _0,r\right) \dot M
    (\vartheta _0,r)}{\varepsilon \sigma \left(r\right)^2}{\rm d}r
  \Bigr|\geq \frac{c_1 \varepsilon^{-\frac{1}{12}}}{4} \right)\\
&\qquad \qquad\qquad \qquad +\Pb_{\vartheta
    _0}\left(\frac{\left|\check\vartheta _{\tau _\varepsilon ,\varepsilon
    }-\vartheta _0\right|^2}{\sqrt{\varepsilon } }\geq \frac{c \varepsilon^{-\frac{1}{12}}}{4}\right)\\
 &\qquad \leq
  \Pb_{\vartheta _0}\left(\frac{\left|\check\vartheta _{\tau _\varepsilon
      ,\varepsilon }-\vartheta _0\right|}{\sqrt{\varepsilon } }\sup_{t-\delta
    \leq s\leq t}\Bigl|\int_{\tau _\varepsilon }^{s}  Q\left(\vartheta _0,r\right) \left(\xi
_{r,\varepsilon }^2-\frac{1}{2}\right)   {\rm d}r
  \Bigr|\geq \frac{c_1 \varepsilon^{-\frac{1}{12}}}{2} \right)\\
\end{align*}
where
\begin{align*}
Q\left(\vartheta _0,r\right)= \frac{ S(\vartheta _0 ,r)
}{ \sigma \left(r\right)}{ }\left[\frac{\dot f(\vartheta _0  ,r)}{f(\vartheta _0 
    ,r)}+\frac{\dot b(\vartheta _0 
    ,r)}{b(\vartheta _0 ,r)} \right]^2.
\end{align*}
We have 
\begin{align*}
&\sup_{t-\delta
    \leq s\leq t}\Bigl|\int_{\tau _\varepsilon }^{s}  Q\left(\vartheta _0,r\right) \left(\xi
_{r,\varepsilon }^2-\frac{1}{2}\right)   {\rm d}r
  \Bigr|\\
&\qquad \qquad\leq\sup_{t-\delta
    \leq s\leq t}\Bigl|\int_{s}^{t}  Q\left(\vartheta _0,r\right) \left(\xi
_{r,\varepsilon }^2-\frac{1}{2}\right)   {\rm d}r
  \Bigr|  +\Bigl|\int_{\tau _\varepsilon }^{t}  Q\left(\vartheta _0,r\right) \left(\xi
_{r,\varepsilon }^2-\frac{1}{2}\right)   {\rm d}r
  \Bigr|\\
&\qquad \qquad\leq \delta \Bigl|  Q\left(\vartheta _0,\tilde r\right) \left(\xi
_{\tilde r,\varepsilon }^2-\frac{1}{2}\right)\Bigr|
    +\Bigl|\int_{\tau _\varepsilon }^{t}  Q\left(\vartheta _0,r\right) \left(\xi
_{r,\varepsilon }^2-\frac{1}{2}\right)   {\rm d}r
  \Bigr|.
\end{align*}
Remind that $\delta =\delta _\varepsilon =\varepsilon ^{\frac{7}{8}}$. This
relations together with \eqref{23-56}, \eqref{23-77} and the estimate (21) in \cite{Kut19}
allow us to write for any $p_*>0$
\begin{align*}
\Pb_{\vartheta _0}\left(\sup_{t-\delta \leq s\leq t}\frac{\left|\vartheta
_{s,\varepsilon }^\star-\vartheta 
_0\right|}{\sqrt{\varepsilon }}\geq \varepsilon^{-\frac{1}{12}}  \right) \leq
C\,\varepsilon ^{p_*} .
\end{align*}

\end{proof}
\begin{remark}
\label{R23-5} {\rm The expression for $\Ex_{\vartheta _0}Y_t^2$ can be
  immediately obtained from the representation
\begin{align*}
Y_t=y_0\exp\left(\int_{0}^{t}a\left(\vartheta _0,v\right){\rm d}v\right)+
\int_{0}^{t}\exp\left(\int_{s}^{t}a\left(\vartheta _0,v\right){\rm
  d}v\right)b\left(\vartheta _0,s\right) {\rm d}V_s.
\end{align*}
}
\end{remark}

\begin{remark}
\label{R23-6} {\rm 

We have the weak convergences
\begin{align}
\label{23-78}
\frac{Y_t-m\left(\vartheta _0,t\right)}{\sqrt{\varepsilon }}\Longrightarrow \sqrt{\frac{ b\left(\vartheta _0,t
    \right)\sigma \left(t\right)}{{f\left(\vartheta _0,t\right) }}}\;\left(\xi
_{t }^{\left(2\right)}-\xi _{t}^{\left(1\right)}\right)
\end{align}
and 
\begin{align}
\label{23-79}
&\frac{Y_t-\hat m_t}{\sqrt{\varepsilon }} \Longrightarrow\frac{\dot f\left(\vartheta
  _0,t\right) }{f\left(\vartheta
  _0,t\right)}\;\hat\eta_{t}(\vartheta _0)\,Y_t+ \sqrt{\frac{ b\left(\vartheta _0,t
    \right)\sigma \left(t\right)}{{f\left(\vartheta _0,t\right) }}}\;\left(\xi
_{t }^{\left(2\right)}-\xi _{t}^{\left(1\right)}\right).
\end{align}
Recall that  the Gaussian process $ \hat\eta_{t}(\vartheta _0),0<t\leq T$ and
Gaussian r.v.'s $\xi _t^{\left(1\right)},\xi _t^{\left(2\right)},0<t\leq T$ are independent.

The variances of the limit laws  are
\begin{align*}
\gamma _*\left(\vartheta _0,t\right)=\frac{ b\left(\vartheta _0,t
    \right)\sigma \left(t\right)}{{f\left(\vartheta _0,t\right) }}
\end{align*}
in the first case and 
\begin{align*}
\Gamma_* \left(\vartheta _0,t\right)={\frac{ b\left(\vartheta _0,t
    \right)\sigma \left(t\right)}{{f\left(\vartheta _0,t\right) }}}+\frac{\dot f\left(\vartheta
  _0,t\right)^2 \Ex_{\vartheta _0}Y_t^2}{f\left(\vartheta
  _0,t\right)^2{\rm I}_0^t(\vartheta _0)}=\gamma _*\left(\vartheta _0,t\right)+\frac{\dot f\left(\vartheta
  _0,t\right)^2\Ex_{\vartheta _0}Y_t^2 }{f\left(\vartheta
  _0,t\right)^2{\rm I}_0^t(\vartheta _0)}
\end{align*}
in the second case.

It is interesting to note that if $f\left(\vartheta ,t\right)=f\left(t\right)
$ and the function $b\left(\vartheta ,t\right)$ depends on unknown parameter
$\vartheta $ then  the limit representation and the limit variance in both cases
coincide. This means that {\it the  unobserved component $Y_t$ is approximated by
$\hat m_t$ with the  same asymptotic precision as if the value $\vartheta $ is
known}. 

For the integral linear and quadratic  functionals we have the following
limits : for any $t_0>0$
\begin{align*}
&\int_{t_0}^{T}\frac{Y_t-m\left(\vartheta _0,t\right)}{\sqrt{\varepsilon
}}\,{\rm d}t \longrightarrow 0,\\
&\int_{t_0}^{T}\left(\frac{Y_t-m\left(\vartheta _0,t\right)}{\sqrt{\varepsilon
}}\right)^2\,{\rm d}t \longrightarrow \int_{t_0}^{T}\frac{b\left(\vartheta _0,t
    \right)\sigma \left(t\right)}{f\left(\vartheta _0,t\right)}{\rm d}t=\int_{t_0}^{T}\gamma _*\left(\vartheta _0,t\right){\rm d}t,
\end{align*}
and
\begin{align*}
&\int_{t_0}^{T}\frac{Y_t-\hat m_t}{\sqrt{\varepsilon
}}\,{\rm d}t \Longrightarrow  \int_{t_0}^{T}\frac{\dot f\left(\vartheta
  _0,t\right)}{f\left(\vartheta
  _0,t\right)}\;\hat\eta_{t}(\vartheta _0)\,Y_t{\rm d}t   ,\\
&\int_{t_0}^{T}\left(\frac{Y_t-\hat m_t}{\sqrt{\varepsilon
}}\right)^2\,{\rm d}t \Longrightarrow \int_{t_0}^{T} \left[   \frac{\dot f\left(\vartheta
  _0,t\right)^2}{f\left(\vartheta
  _0,t\right)^2}\;\hat\eta_{t}(\vartheta _0)^2\,Y_t^2+ \gamma _*\left(\vartheta _0,t\right) \right]    {\rm d}t.
\end{align*}
The proofs follow from Lemma 4 in \cite{Kut19}.

}
\end{remark}

\begin{remark}
\label{R23-7} {\rm 
Consider the question of asymptotic optimality of the estimation of the
conditional expectation $m\left(\vartheta ,t\right)$. Remind that
$m\left(\vartheta ,t\right)$ itself   is mean squared  optimal
estimator of $Y_t$. 

Let us denote
\begin{align*}
D\left(\vartheta _0,t\right)^2=\frac{\dot f\left(\vartheta _0,t\right)^2}{f\left(\vartheta _0,t\right)^2}\Ex_{\vartheta _0} \left[\hat \eta _t \left(\vartheta
_0\right) Y_t\right] ^2=\frac{\dot f\left(\vartheta
  _0,t\right)^2\,\Ex_{\vartheta _0}Y_t ^2}{f\left(\vartheta _0,t\right)^2{\rm I}_0^t\left(\vartheta
  _0\right)    } .
\end{align*}
We have the lower minimax bound: for any $t\in (0,T]$ and any $\vartheta _0\in\Theta $
\begin{align}
\label{23-80}
\lim_{\nu \rightarrow }\Liminf_{\varepsilon \rightarrow
  0}\sup_{\left|\vartheta -\vartheta _0\right|\leq \nu } \varepsilon
^{-1}\Ex_{\vartheta}\left| \bar m\left(t\right)-m\left(\vartheta
,t\right)\right|^2\geq  D\left(\vartheta _0,t\right)^2,
\end{align}
where $ \bar m\left(t\right) $ is an arbitrary estimator of $m\left(\vartheta
,t\right)$. We follow the proof of Theorem 1.9.1 in \cite{IH81}.  We have
elementary estimates 
\begin{align*}
\sup_{\left|\vartheta -\vartheta _0\right|\leq \nu } \Ex_{\vartheta}\left|
\bar m\left(t\right)-m\left(\vartheta ,t\right)\right|^2&\geq \int_{\vartheta
  _0-\nu }^{\vartheta _0+\nu}\Ex_{\vartheta}\left| \bar
m\left(t\right)-m\left(\vartheta ,t\right)\right|^2p\left(\vartheta
\right){\rm d}\vartheta\\
&\geq \int_{\vartheta
  _0-\nu }^{\vartheta _0+\nu}\Ex_{\vartheta}\left| \tilde
m\left(t\right)-m\left(\vartheta ,t\right)\right|^2p\left(\vartheta
\right){\rm d}\vartheta.
\end{align*}
 Here we introduced a continuous positive density $p\left(\vartheta
 \right),\vartheta \in \left(\vartheta _0-\delta ,\vartheta _0+\delta\right) $
 and $\tilde 
m\left(t\right)$ is Bayesian estimator, which corresponds to this density
$p\left(\cdot \right)$ and quadratic loss function, i.e.,
\begin{align*}
\tilde m\left(t\right)=\int_{\vartheta _0-\nu }^{\vartheta _0+\nu
}m\left(\vartheta ,t\right) p\left(\vartheta |X^t\right){\rm
  d}\vartheta=\frac{\int_{\vartheta _0-\nu }^{\vartheta _0+\nu
  }m\left(\vartheta ,t\right) p\left(\vartheta \right)L\left(\vartheta
  ,X^t\right){\rm d}\vartheta}{\int_{\vartheta _0-\nu }^{\vartheta _0+\nu
  } p\left(\vartheta\right)L\left(\vartheta
  ,X^t\right){\rm d}\vartheta}.
\end{align*}
Below we change the variables  $\vartheta =\vartheta _0+\sqrt{\varepsilon }u$ and
 denote 
$$Z_{t,\varepsilon} \left(u\right)=\frac{L\left(\vartheta _0+\sqrt{\varepsilon }u
  ,X^t\right)}{L\left(\vartheta_0
  ,X^t\right)}   ,\qquad \UU_\varepsilon =\left(-\nu \varepsilon ^{-1/2},\nu
 \varepsilon ^{-1/2}\right).
$$
We can write
\begin{align*}
&\frac{\int_{\vartheta _0-\nu }^{\vartheta _0+\nu }m\left(\vartheta ,t\right)
    p\left(\vartheta \right)L\left(\vartheta ,X^t\right){\rm
      d}\vartheta}{\int_{\vartheta _0-\nu }^{\vartheta _0+\nu }
    p\left(\vartheta\right)L\left(\vartheta ,X^t\right){\rm
      d}\vartheta}=\frac{\int_{\vartheta _0-\nu }^{\vartheta _0+\nu
    }m\left(\vartheta ,t\right) p\left(\vartheta \right)\frac{L\left(\vartheta
      ,X^t\right)}{L\left(\vartheta_0 ,X^t\right) }{\rm
      d}\vartheta}{\int_{\vartheta _0-\nu }^{\vartheta _0+\nu }
    p\left(\vartheta\right)\frac{L\left(\vartheta
      ,X^t\right)}{L\left(\vartheta_0 ,X^t\right) }{\rm d}\vartheta}\\
 &\qquad
  =\int_{\UU_\varepsilon   }^{ } m\left(\vartheta
  _0+\sqrt{\varepsilon }u ,t\right) p\left(\vartheta _0+\sqrt{\varepsilon }u
  \right)Z_{t,\varepsilon}\left(u\right){\rm
    d}u\left(\int_{\UU_\varepsilon  }^{ }
  p\left(\vartheta _0+\sqrt{\varepsilon }u \right)Z_{t,\varepsilon}
  \left(u\right){\rm d}u \right)^{-1}\\
 &\qquad
  = m\left(\vartheta
  _0 ,t\right)+\sqrt{\varepsilon } \dot  m\left(\vartheta
  _0 ,t\right)  \frac{\int_{\UU_\varepsilon }^{ }u
  p\left(\vartheta _0
  \right)Z_{t,\varepsilon} \left(u\right){\rm
    d}u}{ \int_{\UU_\varepsilon  }^{}
  p\left(\vartheta _0\right)Z_{t,\varepsilon}
  \left(u\right){\rm d}u} \left(1+o\left(1\right)\right) \\
 &\qquad
  = m\left(\vartheta
  _0 ,t\right)-\sqrt{\varepsilon } \frac{\dot f\left(\vartheta _0,t\right)}{f\left(\vartheta _0,t\right)}  m\left(\vartheta
  _0 ,t\right)  \frac{\int_{\UU_\varepsilon }^{ }u
 Z_{t,\varepsilon}\left(u\right){\rm
    d}u}{ \int_{\UU_\varepsilon  }^{}
  Z_{t,\varepsilon}
  \left(u\right){\rm d}u } \left(1+o\left(1\right) \right),
\end{align*}
where we  used the expansion 
\begin{align*}
m\left(\vartheta _0+\sqrt{\varepsilon }u ,t\right)&=m\left(\vartheta _0
,t\right)+u \dot m\left(\vartheta_0
,t\right)\,\sqrt{\varepsilon }\left(1+o\left(1\right)\right)\\
&=m\left(\vartheta _0
,t\right)- u\, \frac{\dot f\left(\vartheta _0,t\right)}{f\left(\vartheta
  _0,t\right)}  m\left(\vartheta   _0 ,t\right) \,\sqrt{\varepsilon }\left(1+o\left(1\right)\right).
\end{align*}
Therefore
\begin{align*}
\frac{\tilde m\left(t\right)-m\left(\vartheta _0
,t\right)  }{\sqrt{\varepsilon }}&=-  \frac{\dot f\left(\vartheta _0,t\right)}{f\left(\vartheta
  _0,t\right)}\,  m\left(\vartheta   _0 ,t\right)\, \frac{\int_{\UU_\varepsilon }^{ }u
 Z_{t,\varepsilon}\left(u\right){\rm
    d}u}{ \int_{\UU_\varepsilon  }^{}
  Z_{t,\varepsilon}
  \left(u\right){\rm d}u } \left(1+o\left(1\right) \right)\\
&=-  \frac{\dot f\left(\vartheta _0,t\right)}{f\left(\vartheta
  _0,t\right)}\,\, Y_t\, \frac{\int_{\UU_\varepsilon }^{ }u
 Z_{t,\varepsilon}\left(u\right){\rm
    d}u}{ \int_{\UU_\varepsilon  }^{}
  Z_{t,\varepsilon}
  \left(u\right){\rm d}u } \left(1+o\left(1\right) \right).
\end{align*}
 
In the work \cite{Kut19}, Theorem 1 it was shown that Bayes estimators constructed by
observations $X^T$  are asymptotically
normal 
\begin{align*}
\frac{\tilde\vartheta _\varepsilon -\vartheta _0}{\sqrt{\varepsilon
}}=\frac{\int_{\UU_\varepsilon }^{ }u Z_{T,\varepsilon}\left(u\right){\rm
    d}u}{ \int_{\UU_\varepsilon }^{} Z_{T,\varepsilon} \left(u\right){\rm d}u
} \left(1+o\left(1\right) \right)\Longrightarrow \hat \eta _T\left(\vartheta _0\right)
\end{align*}
and all polynomial moments converge. 

Hence
\begin{align*}
\frac{\tilde m\left(t\right)-m\left(\vartheta _0
,t\right)  }{\sqrt{\varepsilon }}\Longrightarrow -  \frac{\dot f\left(\vartheta _0,t\right)}{f\left(\vartheta
  _0,t\right)}  \,\hat \eta _t\left(\vartheta _0\right)\,Y_t
\end{align*}
and we have the convergence ($\varepsilon \rightarrow 0$)
\begin{align*}
\int_{\vartheta
  _0-\nu }^{\vartheta _0+\nu}\Ex_{\vartheta}\frac{\left| \tilde
m\left(t\right)-m\left(\vartheta ,t\right)\right|^2}{\varepsilon }p\left(\vartheta
\right){\rm d}\vartheta&\longrightarrow \int_{\vartheta
  _0-\nu }^{\vartheta _0+\nu}\frac{\dot f\left(\vartheta,t\right)^2}{f\left(\vartheta
  ,t\right)^2} \Ex_{\vartheta}\Bigr(\hat \eta
_t\left(\vartheta\right)Y_t \Bigl)^2p\left(\vartheta
\right){\rm d}\vartheta \\
&=\int_{\vartheta
  _0-\nu }^{\vartheta _0+\nu}D\left(\vartheta ,t\right)^2p\left(\vartheta
\right){\rm d}\vartheta.
\end{align*}

Therefore we obtained the relation
\begin{align*}
\Liminf_{\varepsilon \rightarrow
  0}\sup_{\left|\vartheta -\vartheta _0\right|\leq \nu } \varepsilon
^{-1}\Ex_{\vartheta}\left| \bar m\left(t\right)-m\left(\vartheta
,t\right)\right|^2\geq \int_{\vartheta
  _0-\nu }^{\vartheta _0+\nu}D\left(\vartheta ,t\right)^2p\left(\vartheta
\right){\rm d}\vartheta.
\end{align*}
The function $D\left(\vartheta ,t\right)^2 $ is continuous on $\vartheta $ and
we have the limit
\begin{align*}
\lim_{\nu \rightarrow 0}\int_{\vartheta
  _0-\nu }^{\vartheta _0+\nu}D\left(\vartheta ,t\right)^2p\left(\vartheta
\right){\rm d}\vartheta=D\left(\vartheta_0 ,t\right)^2. 
\end{align*}
The  bound \eqref{23-80} is proved.

From  relations \eqref{23-78},\eqref{23-79} we obtain
\begin{align}
\label{23-81}
\frac{\hat m_t-m\left(\vartheta_0 ,t\right)}{\sqrt{\varepsilon }}\Longrightarrow -\frac{\dot
  f\left(\vartheta _0,t\right)}{f\left(\vartheta _0,t\right)}\;\hat \eta
_t\left(\vartheta _0\right) \,Y_t.
\end{align}
Therefore the  variance of the limit law is $D\left(\vartheta_0 ,t\right)^2
$. It can be shown that the convergence \eqref{23-81}  is uniform on compacts
  $ \KK\subset\Theta  $ and we have 
\begin{align*}
\lim_{\nu \rightarrow 0}\lim_{\varepsilon \rightarrow
  0}\sup_{\left|\vartheta -\vartheta _0\right|\leq \nu } \varepsilon
^{-1}\Ex_{\vartheta}\left| \hat m_t-m\left(\vartheta
,t\right)\right|^2= D\left(\vartheta_0 ,t\right)^2,
\end{align*}
i.e., $\hat m_t$ is asymptotically efficient (minimax) estimator of
$m\left(\vartheta ,t\right)$. 

}

\end{remark}

{\it Example 3.} Consider the model of observations of Example 2
\begin{align*}
{\rm d}X_t&=f\left(t\right)Y_t{\rm d}t+\varepsilon \sigma\left(t\right){\rm d}W_t,\qquad
X_0=0, \qquad 0\leq t\leq T,\\
{\rm d}Y_t&=a\left(t\right)Y_t{\rm
  d}t+\sqrt{h\left(t\right)+\vartheta g\left(t\right)}{\rm d}V_t, \quad Y_0=y_0
\end{align*}
where the functions $f\left(\cdot \right),\sigma \left(\cdot
\right),a\left(\cdot \right),h\left(\cdot \right),g\left(\cdot \right) $ have
continuous derivatives on $t$ and are strictly positive  and $\vartheta \in
\left(\alpha ,\beta 
\right)$, where $\alpha >0$. It is easy to see that  the conditions
${\scr A},{\scr B},\tilde{\scr C}_1,{\scr C}_2$ in this case are
fulfilled.

The preliminary estimator is:
\begin{align}
\label{23-82}
\check\vartheta _{\tau _\varepsilon ,\varepsilon }=\left(\Psi _{\tau _\varepsilon
  ,\varepsilon }-\int_{0}^{\tau _\varepsilon
}f\left(t\right)^2h\left(t\right){\rm d}t\right)\left( \int_{0}^{\tau _\varepsilon
}f\left(t\right)^2g\left(t\right){\rm d}t\right)^{-1} 
\end{align}
where $\tau _\varepsilon =\varepsilon ^{1/12}$, Fisher information
\begin{align*}
{\rm I}_{\tau _\varepsilon }^t\left(\vartheta \right)=\int_{\tau _\varepsilon
}^{t}\frac{f\left(s\right)g\left(s\right)^2}{8\left[h\left(s\right)+\vartheta
    g\left(s\right)\right]^{3/2}\sigma \left(s\right)}{\rm d}s. 
\end{align*}
The random process $M(\check\vartheta _{\tau _\varepsilon ,\varepsilon
},t)$ is solution of the equation
\begin{align}
\label{23-83}
{\rm d}M(\check\vartheta _{\tau _\varepsilon ,\varepsilon
},t)=\frac{f\left(t\right)\left[h\left(t\right)+\check\vartheta _{\tau _\varepsilon
    ,\varepsilon }g\left(t\right)\right]^{1/2}}{\varepsilon \sigma
  \left(t\right)}\left[{\rm
    d}X_t-M(\check\vartheta _{\tau _\varepsilon ,\varepsilon },t){\rm
    d}t\right] 
\end{align}
on the interval $\left[\tau _\varepsilon ,T\right]$.

For the derivative we have the following equation
\begin{align}
\label{23-84}
{\rm d}\dot M(\check\vartheta _{\tau _\varepsilon ,\varepsilon
},t)& =-\frac{f\left(t\right)\left[h\left(t\right)+\check\vartheta _{\tau _\varepsilon
    ,\varepsilon }g\left(t\right)\right]^{1/2}}{\varepsilon \sigma
  \left(t\right)} \dot M(\check\vartheta _{\tau _\varepsilon ,\varepsilon
},t) {\rm d}t\nonumber\\
&\qquad \qquad+ \frac{f\left(t\right)g\left(t\right)}{2\varepsilon \sigma
  \left(t\right)\left[h\left(t\right)+\check\vartheta _{\tau _\varepsilon
    ,\varepsilon }g\left(t\right)\right]^{1/2} }
\left[{\rm
    d}X_t-M(\check\vartheta _{\tau _\varepsilon ,\varepsilon },t){\rm
    d}t\right] .
\end{align}
 The initial values  are 
defined by the relation like \eqref{23-59}. The One-step MLE-process has the
same form
\begin{align}
\label{23-85}
\vartheta _{t,\varepsilon}^\star=\check\vartheta _{\tau _\varepsilon
  ,\varepsilon }+\frac{1}{{\rm I}_{\tau _\varepsilon }^t(\check\vartheta _{\tau _\varepsilon
  ,\varepsilon })}\int_{\tau _\varepsilon }^{t}\frac{\dot M(\check\vartheta _{\tau _\varepsilon ,\varepsilon },s)}{\varepsilon \sigma \left(s\right)^2}\left[{\rm d}X_s-M(\check\vartheta _{\tau _\varepsilon ,\varepsilon },s){\rm d}s\right],\quad \tau _\varepsilon <t\leq T.
\end{align}
The adaptive filtre  is
\begin{align}
\label{23-86}
{\rm d}\hat m_t=\frac{f\left(t\right)\left[{h\left(t\right)+\vartheta
    _{t,\varepsilon}^\star g\left(t\right) }\right]^{1/2}}{\varepsilon  \sigma \left(t\right)}\left[{\rm d}X_t-\hat m_t{\rm d}t\right],\qquad \tau _\varepsilon <t\leq T.
\end{align}
For the initial value $\hat m_{\tau _\varepsilon }$ we write 
\begin{align}
\label{23-87}
\hat m_{\tau _\varepsilon }&=y_0\phi\left(\check\vartheta _{\tau,\varepsilon  },
\tau \right)+\phi\left(\check\vartheta _{\tau,\varepsilon  } ,\tau \right)
\left(X_\tau H\left(\check\vartheta _{\tau,\varepsilon  } 
,\tau \right)-\int_{0}^{\tau }X_sH'\left(\check\vartheta _{\tau,\varepsilon  } ,s\right){\rm d}s \right).
\end{align}
Here
\begin{align*}
\phi\left(\vartheta, s\right)&=\exp\left(-\int_{0}^{s}q_\varepsilon
\left(\vartheta ,v\right){\rm d}v \right) ,\qquad \quad q_\varepsilon
\left(\vartheta ,t\right)=-a\left(\vartheta ,t\right)+ \frac{\gamma
  \left(\vartheta ,t\right)f\left(\vartheta ,t\right)^2}{\varepsilon ^2\sigma
  \left(t\right)^2},\\ 
 H\left(\vartheta
,s\right)&=\phi\left(\vartheta, s\right)^{-1} \frac{\gamma \left(\vartheta
  ,s\right)f\left(\vartheta ,s\right)}{\varepsilon ^2\sigma \left(s\right)^2}.
\end{align*}

Recall that $\hat m_t $ is asymptotically efficient estimator of the random
process $m\left(\vartheta _0,t\right)$.

\section{Discussion}

Note that the case of multidimensional parameter $\vartheta $ can be
considered following the same steps as it was done in this work. 

The proposed algorithm of adaptive filtration, is numerically
essentially simpler, than the algorithm based on the substitution of the MLE
$\hat\vartheta _{t,\vartheta }0<t\leq T $ because the calculation of
$\hat\vartheta _{t,\vartheta } $ for each $t$ requires solutions of many
equations \eqref{23-57}-\eqref{23-58}, which has to be used for solution of
the maximum likelihood equations. The solution of Riccati equation is needed
for the calculation of the initial values  only and for  one value of
$\vartheta = \check\vartheta _{\tau,\varepsilon  } $.  Another advantage is the recursive structure
of all equations. Of course, One-step MLE-process can be presented in a
recursive form too.

\section*{Acknowledgments} I would like to thank P.Chigansky for usefull
comments.  This research was financially supported by the  by RSF
research project No 20-61-47043.

\end{document}